\documentclass{article}
\usepackage{amsmath}
\usepackage{amssymb}

\usepackage{theorem}
\numberwithin{equation}{section}

\newtheorem{theorem}{Theorem}[section]
\newtheorem{conjecture}{Conjecture}[section]

\newtheorem{lemma}[theorem]{Lemma}
\newtheorem{definition}{Definition}[section]

\begin{document}

\title{Global well - posedness and scattering for the focusing, energy - critical nonlinear Schr{\"o}dinger problem in dimension $d = 4$ for initial data below a ground state threshold}

\author{Benjamin Dodson}

\maketitle

\noindent \textbf{Abstract:} In this paper we prove global well - posedness and scattering for the focusing, energy - critical nonlinear Schr{\"o}dinger initial value problem in four dimensions. Previous work proved this in five dimensions and higher using the double Duhamel trick. In this paper, using long time Strichartz estimates we are able to overcome the logarithmic blowup in four dimensions.

\section{Introduction}
In this paper we study the Schr{\"o}dinger initial value problem

\begin{equation}\label{equation}
\aligned
i u_{t} + \Delta u &= F(u) = - |u|^{2} u, \\
u(0,x) &= u_{0} \in \dot{H}^{1}(\mathbf{R}^{4}).
\endaligned
\end{equation}

\noindent $(\ref{equation})$ belongs to a class of problems known as the focusing, nonlinear Schr{\"o}dinger initial value problems,

\begin{equation}\label{1.2}
\aligned
i u_{t} + \Delta u &= F(u) = - |u|^{p} u, \\
u(0,x) &= u_{0} \in \dot{H}^{1}(\mathbf{R}^{d}).
\endaligned
\end{equation}

\noindent $(\ref{1.2})$ is called energy - critical if $p = \frac{4}{d - 2}$, $d \geq 3$. In general a solution to $(\ref{1.2})$ conserves the quantities mass,

\begin{equation}\label{1.3}
M(u(t)) = \int |u(t,x)|^{2} dx = M(u(0)),
\end{equation}

\noindent and energy,

\begin{equation}\label{1.4}
E(u(t)) = \frac{1}{2} \int |\nabla u(t,x)|^{2} dx - \frac{1}{p + 2} \int |u(t,x)|^{p + 2} dx = E(u(0)).
\end{equation}

\noindent $(\ref{1.2})$ is called energy - critical when $p = \frac{4}{d - 2}$ since a solution to $(\ref{1.2})$ is invariant under the scaling

\begin{equation}\label{1.5}
u(t,x) \mapsto \lambda^{\frac{d}{2} - 1} u(\lambda^{2} t, \lambda x),
\end{equation}

\noindent and $(\ref{1.5})$ preserves the energy $(\ref{1.4})$.\vspace{5mm}

\noindent There also exist the defocusing, energy - critical problems ($F(u) = |u|^{\frac{4}{d - 2}} u$), which are similar to the focusing problem in some ways, but also contain many important differences. The defocusing problem is now completely worked out.

\begin{theorem}\label{t1.1}
The defocusing initial value problem $(\ref{1.2})$, $F(u) = |u|^{\frac{4}{d - 2}} u$, is globally well - posed and scattering for all $u_{0} \in \dot{H}^{1}(\mathbf{R}^{d})$, $d \geq 3$.
\end{theorem}

\begin{definition}[Scattering]\label{d1.5}
A solution $u$ to $(\ref{1.2})$, $p = \frac{4}{d - 2}$ is said to scatter forward in time if there exists $u_{+} \in \dot{H}^{1}$ such that

\begin{equation}\label{1.16}
\lim_{t \nearrow +\infty} \| u(t) - e^{it \Delta} u_{+} \|_{\dot{H}^{1}(\mathbf{R}^{d})} = 0.
\end{equation}

\noindent Likewise, $u$ is said to scatter backward in time if there exists $u_{-} \in \dot{H}^{1}$ such that

\begin{equation}\label{1.17}
\lim_{t \searrow -\infty} \| u(t) - e^{it \Delta} u_{-} \|_{\dot{H}^{1}(\mathbf{R}^{d})} = 0.
\end{equation}
\end{definition}

\noindent \emph{Proof:} The proof of theorem $\ref{t1.1}$ has involved contributions from a variety of authors. \cite{CaWe1} proved theorem $\ref{t1.1}$ for small data in both the focusing and defocusing problem. \cite{CaWe1} also proved that $(\ref{1.2})$ has a local solution for any initial data $u_{0} \in \dot{H}^{1}(\mathbf{R}^{d})$, where the time of existence depends on the size and profile of $u_{0}$.\vspace{5mm}

\noindent For large data, the seminal result was the work of \cite{B4} and \cite{B3}, proving theorem $\ref{t1.1}$ for radial data in dimensions $d = 3, 4$, and also that for more regular $u_{0}$, this additional smoothness is preserved. See \cite{Gril} for another proof of this last fact. \cite{TerryTao} then extended theorem $\ref{t1.1}$ to radial data in higher dimensions.\vspace{5mm}

\noindent Then \cite{CKSTT4} extended theorem $\ref{t1.1}$ to general $u_{0} \in \dot{H}^{1}$ when $d = 3$. Subsequently, \cite{RhV} extended this to dimension $d = 4$, and \cite{Visan}, \cite{V2} extended theorem $\ref{t1.1}$ to dimensions $d \geq 5$. $\Box$\vspace{5mm}

\noindent \textbf{Remark:} \cite{Visan1} and \cite{KV2} reproved theorem $\ref{t1.1}$ in dimensions three and four using the long time Strichartz estimates of \cite{D2}. We will use long time Strichartz estimates similar to the estimates of \cite{KV2} in this paper as well.\vspace{5mm}

\noindent Returning to the focusing problem, we remark that theorem $\ref{t1.1}$ does not hold for arbitrary data. In fact, by the virial identity (see for example \cite{Glassey})

\begin{equation}\label{1.6}
\frac{d^{2}}{dt^{2}} \int |x|^{2} |u(t,x)|^{2} dx = 8 [\int |\nabla u(t,x)|^{2} dx - \int |u(t,x)|^{\frac{2d}{d - 2}} dx],
\end{equation}

\noindent so for $x u_{0} \in L^{2}(\mathbf{R}^{d})$ and $E(u_{0}) < 0$, the solution must break down in finite time. Moreover,

\begin{equation}\label{1.7}
W(x) = W(x,t) = \frac{1}{(1 + \frac{|x|^{2}}{d(d - 2)})^{\frac{d - 2}{2}}}
\end{equation}

\noindent lies in $\dot{H}^{1}(\mathbf{R}^{d})$ and solves the elliptic equation

\begin{equation}\label{1.8}
\Delta W + |W|^{\frac{4}{d - 2}} W = 0.
\end{equation}

\noindent Therefore, scattering cannot always occur even for global solutions. Instead, as in the mass - critical problem, we conjecture that scattering holds for initial data below the threshold given by $(\ref{1.7})$.\vspace{5mm}

\begin{conjecture}\label{c1.2}
Let $d \geq 3$ and let $u : I \times \mathbf{R}^{d} \rightarrow \mathbf{C}$ be a solution to $(\ref{1.2})$, $p = \frac{4}{d - 2}$. If

\begin{equation}\label{1.9}
\| u_{0} \|_{\dot{H}^{1}(\mathbf{R}^{d})} < \| W \|_{\dot{H}^{1}(\mathbf{R}^{d})},
\end{equation}

\noindent and 

\begin{equation}\label{1.10}
E(u_{0}) < E(W),
\end{equation}

\noindent then

\begin{equation}\label{1.11}
\int_{I} \int |u(t,x)|^{\frac{2(d + 2)}{d - 2}} dx dt \leq C(\| u_{0} \|_{\dot{H}^{1}}, E(u_{0})) < \infty.
\end{equation}
\end{conjecture}

\noindent \cite{CaWe} and \cite{CaWe1} proved that $(\ref{1.2})$, $p = \frac{4}{d - 2}$ is well - posed on $I$ for initial data $u_{0}$ if and only if, for any $J \subset I$ compact, $S_{J}(u) < \infty$. If $S_{[t_{1}, \infty)}(u) < \infty$ for some $t_{1} \in \mathbf{R}$, then $u$ scatters forward in time. Likewise, if $S_{(-\infty, t_{1}]}(u) < \infty$ then $u$ scatters backward in time.

\begin{definition}[Scattering size]\label{d1.3}
The scattering size of a solution to $(\ref{1.2})$ on a time interval $I$ is given by

\begin{equation}\label{1.13}
S_{I}(u) = \int_{I} \int_{\mathbf{R}^{d}} |u(t,x)|^{\frac{2(d + 2)}{d - 2}} dx dt.
\end{equation}
\end{definition}

\begin{definition}[Blow up]\label{d1.4}
A solution $u$ to $(\ref{1.2})$ blows up forward in time on $I$ if there exists $t_{1} \in I$ such that

\begin{equation}\label{1.14}
S_{[t_{1}, \sup(I))}(u) = \infty.
\end{equation}

\noindent $u$ blows up backward in time if there exists $t_{1} \in I$ such that

\begin{equation}\label{1.15}
S_{(\inf(I), t_{1}]}(u) = \infty.
\end{equation}
\end{definition}

\noindent Substantial progress has been made toward the proof of conjecture $\ref{c1.2}$.

\begin{theorem}\label{t1.6}
Assume that $E(u_{0}) < E(W)$, $\| u_{0} \|_{\dot{H}^{1}} < \| W \|_{\dot{H}^{1}}$, $d = 3, 4, 5$, and $u_{0}$ is radial. Then $(\ref{1.2})$ is globally well - posed and scatters forward and backward in time.
\end{theorem}

\noindent \emph{Proof:} See \cite{KM1}. $\Box$\vspace{5mm}

\noindent Then \cite{KV1} treated the nonradial case.

\begin{theorem}\label{t1.7}
Assume that $E(u_{0}) < E(W)$, $\| u_{0} \|_{\dot{H}^{1}} < \| W \|_{\dot{H}^{1}}$, $d \geq 5$. Then $(\ref{1.2})$ is globally well - posed and scatters forward and backward in time.
\end{theorem}

\noindent \emph{Proof:} See \cite{KV1}. $\Box$\vspace{5mm}

\noindent \textbf{Remark:} The result of \cite{KV1} was proved under the assumption that

\begin{equation}\label{1.19}
\| u \|_{L_{t}^{\infty} \dot{H}_{x}^{1}(I \times \mathbf{R}^{4})} < \| \nabla W \|_{L^{2}(\mathbf{R}^{4})}.
\end{equation}

\noindent Now by the energy trapping lemma of \cite{KM1}, if $E(u_{0}) < E(W)$ and $\| u_{0} \|_{\dot{H}^{1}} < \| W \|_{\dot{H}^{1}}$, $(\ref{1.19})$ holds.\vspace{5mm}

\begin{lemma}\label{l1.12}
If $E(u_{0}) \leq (1 - \delta) E(W)$ and $\| \nabla u_{0} \|_{L^{2}(\mathbf{R}^{d})} < (1 - \delta) \| \nabla W \|_{L^{2}(\mathbf{R}^{d})}$ for some $\delta > 0$, then there exists $\bar{\delta}(\delta, d) > 0$ such that for all $t \in I$, where $I$ is the maximal interval of existence of $u$,

\begin{equation}\label{1.20}
\| \nabla u(t) \|_{L_{x}^{2}(\mathbf{R}^{d})} \leq (1 - \bar{\delta}) \| \nabla W \|_{L^{2}(\mathbf{R}^{d})}.
\end{equation}
\end{lemma}

\noindent \emph{Proof:} This follows from the work of \cite{Aubin} and \cite{Talenti}, which proved that if $C_{d}$ is the best constant in the Sobolev embedding,

\begin{equation}\label{1.22}
\| u \|_{L_{x}^{\frac{2d}{d - 2}}(\mathbf{R}^{d})} \leq C_{d} \| \nabla u \|_{L_{x}^{2}(\mathbf{R}^{d})},
\end{equation}

\noindent and

\begin{equation}\label{1.23}
\| u \|_{L_{x}^{\frac{2d}{d - 2}}(\mathbf{R}^{d})} = C_{d} \| \nabla u \|_{L_{x}^{2}(\mathbf{R}^{d})},
\end{equation}

\noindent then $u = C W_{\theta_{0}, x_{0}, \lambda_{0}}$ for some constant $C$, $\theta_{0} \in \mathbf{R}$, $x_{0} \in \mathbf{R}^{d}$, and $\lambda_{0} \in (0, \infty)$, and

\begin{equation}\label{1.24}
W_{\theta_{0}, x_{0}, \lambda_{0}} = \frac{1}{\lambda_{0}^{\frac{d - 2}{2}}} e^{i \theta_{0}} W(\frac{x - x_{0}}{\lambda_{0}}),
\end{equation}

\noindent $W$ is given by $(\ref{1.7})$. In particular, when $d = 4$, $(\ref{1.8})$ implies

\begin{equation}\label{1.25}
0 = \langle \Delta W, W \rangle + \langle W, |W|^{2} W \rangle = -\int |\nabla W|^{2} dx + \int |W|^{4} dx.
\end{equation}

\noindent Then by $(\ref{1.23})$,

\begin{equation}\label{1.26}
C_{4} = \frac{1}{\| W \|_{L_{x}^{4}(\mathbf{R}^{4})}},
\end{equation}

\noindent so

\begin{equation}\label{1.27}
E(W)(1 - \delta) \geq E(u_{0}) = \frac{1}{2} \int |\nabla u(t)|^{2} dx (1 - \frac{1}{2} \frac{\| u(t) \|_{L_{x}^{4}(\mathbf{R}^{4})}^{2}}{\| W \|_{L_{x}^{4}(\mathbf{R}^{4})}^{2}}).
\end{equation}

\noindent Now make a bootstrap argument. Since $\| u(0) \|_{\dot{H}^{1}} < \| W \|_{\dot{H}^{1}}$, by local well - posedness 

\begin{equation}\label{1.28}
\| u(t) \|_{\dot{H}^{1}} \leq \| \nabla W \|_{L^{2}}
\end{equation}

\noindent on some closed interval $J$ of $0$. Then by $(\ref{1.23})$, $(\ref{1.26})$, and $(\ref{1.27})$,

\begin{equation}\label{1.29}
E(W) (1 - \delta) = (1 - \delta) \frac{1}{4} \| W \|_{\dot{H}^{1}(\mathbf{R}^{4})}^{2} \geq \frac{1}{4} \| \nabla u(t) \|_{L_{x}^{2}(\mathbf{R}^{4})}^{2},
\end{equation}

\noindent which in turn implies that $\| u(t) \|_{\dot{H}^{1}(\mathbf{R}^{4})}^{2} \leq (1 - \delta) \| W \|_{\dot{H}^{1}(\mathbf{R}^{4})}^{2}$. $\Box$\vspace{5mm}

\noindent Scattering results for the mass - critical problem (\cite{Merle}, \cite{KVZ}, \cite{TVZ}, \cite{D5}) assume that the initial data $u_{0}$ has mass below the mass of a ground state. For the energy - critical problem it stands to reason that there should be two assumptions on the initial data because unlike the mass $(\ref{1.3})$ , the $\dot{H}^{1}$ norm is not conserved. On the other hand, while energy is conserved, energy is not positive definite $(\ref{1.4})$, so $E(u(t)) < E(W)$ does not by itself give a bound on the size of $u(t)$. The author of this paper is personally unaware of any solutions $u(t)$ to $(\ref{1.2})$, $p = \frac{4}{d - 2}$ that satisfy $(\ref{1.19})$ but not the initial conditions of theorem $\ref{t1.6}$, although he suspects that there most likely are.\vspace{5mm}

\noindent In this paper we prove global well - posedness and scattering for nonradial data in dimension four.

\begin{theorem}\label{t1.8}
Assume that $E(u_{0}) < E(W)$, $\| u_{0} \|_{\dot{H}^{1}} < \| W \|_{\dot{H}^{1}}$, and $d = 4$. Then $(\ref{1.2})$ is globally well - posed and scatters forward and backward in time.
\end{theorem}

\noindent As in \cite{KM1} and \cite{KV1}, the proof uses the concentration compactness method.

\begin{theorem}\label{t1.9}
If $(\ref{equation})$ is not globally well - posed and scattering for all data satisfying $\| u_{0} \|_{\dot{H}^{1}} < \| W \|_{\dot{H}^{1}}$ and $E(u_{0}) < E(W)$, then there exists a nonzero solution $u$ to $(\ref{equation})$ on $I$, where $I$ is the maximal interval of its existence, such that $u$ is almost periodic for all $t \in I$.
\end{theorem}

\begin{definition}[Almost periodicity]\label{d1.10}
$u(t)$ is said to be almost periodic for all $t \in I$ if there exists $N(t) : I \rightarrow (0, \infty)$ and $x(t) : I \rightarrow \mathbf{R}^{4}$ such that $\frac{1}{N(t)} u(\frac{x - x(t)}{N(t)})$ lies in a compact set $K \subset \dot{H}^{1}(\mathbf{R}^{4})$ for all $t \in I$.
\end{definition}

\begin{theorem}\label{t1.10}
The only almost periodic solution to $(\ref{equation})$ on the maximal interval of its existence $I$, with $\| \nabla u(t) \|_{L_{t}^{\infty} L_{x}^{2}(I \times \mathbf{R}^{4})} < \| \nabla W \|_{L^{2}}$, is $u \equiv 0$.
\end{theorem}

\noindent Then to prove theorem $\ref{t1.8}$, it suffices to show theorems $\ref{t1.9}$ and $\ref{t1.10}$. In fact,

\begin{theorem}\label{t1.11}
To prove theorem $\ref{t1.10}$ it suffices to show that the only global, almost periodic solution to $(\ref{equation})$ on $\mathbf{R}$ with

\begin{equation}\label{1.18}
N(t) \geq 1, \hspace{5mm} N(0) = 1,
\end{equation}

\noindent is $u \equiv 0$.
\end{theorem}

\noindent The main difference between \cite{KV1} and this result in dimension $d = 4$ is that in dimensions $d \geq 5$ the dispersive estimate $(\ref{2.13})$ is doubly integrable, allowing \cite{KV1} to make use of the double Duhamel trick. However, here, even though we can prove $u \in L_{t}^{\infty} L_{x}^{3}$, and thus $F(u) \in L^{1}$, the double integral of $(\ref{2.14})$ diverges logarithmically. Nevertheless, this logarithmically divergent result is good enough to be used in an interaction Morawetz estimate, proving theorem $\ref{t1.11}$.\vspace{5mm}

\noindent \textbf{Outline of Proof:} In $\S 2$, some linear estimates and harmonic analysis results will be discussed. These results will be used frequently throughout the rest of the paper. Only one of the results in this section is new.\vspace{5mm}

\noindent In $\S 3$, the concentration compactness method will be discussed, sketching \cite{KM1} and then \cite{KV1}'s proof of theorems $\ref{t1.9}$. We will also discuss almost periodic solutions to $(\ref{equation})$ and sketch \cite{KV1}'s proof of $\ref{t1.11}$. Finally we will bound the $L_{t}^{\infty} L_{x}^{3}(\mathbf{R} \times \mathbf{R}^{4})$ norm of a solution satisfying $(\ref{1.18})$. In $\S 4$ we prove the long time Strichartz estimate. In contrast to \cite{Visan1} and \cite{KV3}, we will consider the quantity

\begin{equation}\label{1.19}
\int_{I} \frac{1}{N(t)^{2}} dt.
\end{equation}

\noindent The long time Strichartz estimates allow us to easily exclude the case when $\int_{\mathbf{R}} N(t)^{-2} dt < \infty$. In $\S 5$ we show that the soliton blowup solution, that is $N(t) \equiv 1$, is $u \equiv 0$. Finally, in $\S 6$ we will extend this argument to a quasi soliton solution, $(\ref{1.19}) = \infty$. This completes the proof of theorem $\ref{t1.8}$.\vspace{5mm}

\noindent \textbf{Acknowledgements:} During the time of researching this paper, the author was supported by NSF postdoctoral fellowship DMS - 1103914. The author also performed much of the research while a guest of the Hausdorff Institute at the University of Bonn for the summer trimester program in harmonic analysis and partial differential equations. The author is grateful to Rowan Killip, Jason Murphy, and Monica Visan for several helpful discussions regarding this problem.\vspace{5mm}

\section{Linear Estimates and harmonic analysis}
\noindent In this section we describe the tools from harmonic analysis that will be used in this paper. None of the results of this section, with the exception of theorem $\ref{maximalStrichartz}$, are new. Theorem $\ref{maximalStrichartz}$ was proved by \cite{KV2} for dimension $d = 3$ only.

\begin{definition}[Fourier transform]\label{d2.1}
Suppose $f \in L^{1}(\mathbf{R}^{d})$. Then

\begin{equation}\label{2.1}
\mathcal F f(\xi) = (2 \pi)^{-d/2} \int e^{-ix \cdot \xi} f(x) dx.
\end{equation}

\noindent The inverse Fourier transform is then given by

\begin{equation}\label{2.2}
\mathcal F^{-1} \hat{f}(x) = (2 \pi)^{-d/2} \int e^{ix \cdot \xi} \hat{f}(\xi) d\xi.
\end{equation}
\end{definition}

\noindent Plancherel's theorem proved that the Fourier transform and inverse Fourier transform provide a unitary transformation between functions in $L_{x}^{2}(\mathbf{R}^{d})$ and functions in $L_{\xi}^{2}(\mathbf{R}^{d})$. Because of this fact it is useful to decompose a function via a partition of unity in Fourier space, or a Littlewood - Paley decomposition.

\begin{definition}[Littlewood - Paley decomposition]\label{d2.2}
Let $\phi \in C_{0}^{\infty}(\mathbf{R}^{d})$ be a radial, decreasing function, $\phi(x) = 1$ for $|x| \leq 1$, $\phi(x)$ is supported on $|x| > 2$. Then for any $j \in \mathbf{Z}$ let

\begin{equation}\label{2.3}
P_{j} f = (2 \pi)^{-d/2} \int e^{ix \cdot \xi} [\phi(2^{-j - 1} \xi) - \phi(2^{-j} \xi)] \hat{f}(\xi) d\xi.
\end{equation}
\end{definition}

\noindent \textbf{Remark:} It is often convenient to write $P_{N}$, which is given by the multiplier

\begin{equation}
[\phi(\frac{1}{N} \xi) - \phi(\frac{1}{2N} \xi)],
\end{equation}

\noindent or to sum over $N \geq M$, which in this case would be over $M = 2^{j} N$, $j \geq 0$.\vspace{5mm}

\noindent To simplify notation we often write $u_{k}$ or $u_{N}$ instead of $P_{k} u$ or $P_{N} u$.

\begin{theorem}[Littlewood - Paley theorem]\label{t2.3}
For any $1 < p < \infty$,

\begin{equation}\label{2.4}
\| (\sum_{j} |P_{j} f|^{2})^{1/2} \|_{L_{x}^{p}(\mathbf{R}^{d})} \sim_{p, d} \| f \|_{L^{p}(\mathbf{R}^{d})}.
\end{equation}
\end{theorem}

\noindent \emph{Proof:} This is a well - known fact from harmonic analysis. See \cite{St1}, \cite{St}, \cite{T3}, or many other sources. $\Box$\vspace{5mm}

\noindent The proof of theorem $\ref{t2.3}$ utilizes the maximal function, which can be defined in any dimension. We will use the maximal function in one dimension only.

\begin{definition}[Maximal function]\label{d2.4}
For a function $f \in L^{p}(\mathbf{R})$, $1 \leq p \leq \infty$,

\begin{equation}\label{2.5}
\mathcal M(f)(x) = \sup_{T > 0} \frac{1}{T} \int_{x - T}^{x + T} |f(t)| dt.
\end{equation}
\end{definition}

\begin{theorem}[Maximal theorem]\label{t2.5}
For any $1 < p \leq \infty$,

\begin{equation}\label{2.6}
\| \mathcal M(f) \|_{L^{p}(\mathbf{R})} \lesssim_{p} \| f \|_{L^{p}(\mathbf{R})}.
\end{equation}
\end{theorem}

\noindent \emph{Proof:} See \cite{St1}, \cite{St}, or \cite{T3}. The proof there is described in any dimension. $\Box$\vspace{5mm}

\begin{theorem}[Sobolev embedding]\label{t2.6}
For $1 \leq p \leq q \leq \infty$,

\begin{equation}\label{2.7}
\| P_{j} f \|_{L^{q}(\mathbf{R}^{d})} \lesssim 2^{jd(\frac{1}{p} - \frac{1}{q})} \| P_{j} f \|_{L^{p}(\mathbf{R}^{d})}.
\end{equation}
\end{theorem}

\noindent \emph{Proof:} See for example \cite{T1}. $\Box$

\begin{lemma}[Bernstein's lemma]\label{l2.7}
For any $s \in \mathbf{R}$, $j \in \mathbf{Z}$, $1 < p < \infty$,

\begin{equation}\label{2.8}
\| P_{j} f \|_{L^{p}(\mathbf{R}^{d})} \sim_{p, d} \| |\nabla|^{s} f \|_{L^{p}(\mathbf{R}^{d})}.
\end{equation}
\end{lemma}

\noindent \emph{Proof:} See \cite{T4}. $\Box$\vspace{5mm}

\noindent Theorem $\ref{t2.3}$, theorem $\ref{t2.6}$, and lemma $\ref{l2.7}$ will be used throughout this paper, frequently in conjunction with one another.\vspace{5mm}

\noindent The Fourier transform is extremely useful to the study of the linear Schr{\"o}dinger problem,

\begin{equation}\label{2.9}
(i \partial_{t} + \Delta) u = F, \hspace{5mm} u(0,x) = u_{0},
\end{equation}

\noindent because the solution to $(\ref{2.9})$ when $F = 0$ is given by

\begin{equation}\label{2.10}
e^{it \Delta} u_{0} = (2 \pi)^{-d/2} \int e^{-it |\xi|^{2}} e^{ix \cdot \xi} \hat{f}(\xi) d\xi,
\end{equation}

\noindent and the general strong solution to $(\ref{2.9})$ is given by

\begin{equation}\label{2.11}
u(t) = e^{i(t - t_{0}) \Delta} u(t_{0}) - i \int_{t_{0}}^{t} e^{i(t - \tau) \Delta} F(\tau) d\tau.
\end{equation}

\noindent Since $|e^{it |\xi|^{2}}| = 1$, 

\begin{equation}\label{2.12}
\| e^{it \Delta} f \|_{L^{2}(\mathbf{R}^{d})} = \| f \|_{L^{2}(\mathbf{R}^{d})},
\end{equation}

\noindent and in fact, for any $L^{2}$ - based Sobolev space,

\begin{equation}\label{2.12.1}
\| e^{it \Delta} f \|_{\dot{H}_{x}^{s}(\mathbf{R}^{d})} = \| f \|_{\dot{H}^{s}(\mathbf{R}^{d})}.
\end{equation}

\noindent By completing the square in the exponent of $(\ref{2.10})$ and stationary phase computations,

\begin{equation}\label{2.13}
e^{it \Delta} f(x) = \frac{1}{(4 \pi t)^{d/2}} e^{-i d\pi/4} \int e^{-i \frac{|x - y|^{2}}{4t}} f(y) dy.
\end{equation}

\noindent \textbf{Remark:} More generally, if $P$ is any Fourier multiplier, that is,

\begin{equation}\label{2.13.1}
P f(x) = \mathcal F^{-1} (P \mathcal Ff)(\xi),
\end{equation}

\noindent then

\begin{equation}\label{2.13.2}
Pf(x) = \int (\mathcal F^{-1} P)(x - y) f(y) dy.
\end{equation}

\noindent Therefore,

\begin{equation}\label{2.14}
\| e^{it \Delta} f \|_{L_{x}^{\infty}(\mathbf{R}^{d})} \lesssim_{d} t^{-d/2} \| f \|_{L^{1}(\mathbf{R}^{d})}.
\end{equation}

\noindent Using both analysis on the Fourier side $(\ref{2.10})$ (see \cite{Stri}), and on the spatial side $(\ref{2.13})$ (see \cite{GV}, \cite{KT}, and \cite{Yaj}), we have the sharp result

\begin{theorem}[Strichartz estimates]\label{t2.8}
For $d \geq 3$, and $(p_{1}, q_{1})$, $(p_{2}, q_{2})$ satisfying $p_{j} \geq 2$,

\begin{equation}\label{2.15}
\frac{2}{p_{j}} = d(\frac{1}{2} - \frac{1}{q_{j}}),
\end{equation}

\noindent if $u$ solves $(\ref{2.9})$ on $I$, $t_{0} \in I$, and $\frac{1}{p'} = 1 - \frac{1}{p}$, then

\begin{equation}\label{2.16}
\| u \|_{L_{t}^{p_{1}} L_{x}^{q_{1}}(I \times \mathbf{R}^{d})} \lesssim_{d} \| u(t_{0}) \|_{L_{x}^{2}(\mathbf{R}^{d})} + \| F \|_{L_{t}^{p_{2}'} L_{x}^{q_{2}'}(I \times \mathbf{R}^{d})}.
\end{equation}
\end{theorem}

\noindent \emph{Proof:} See \cite{Stri} for the seminal result, \cite{GV} and \cite{Yaj} for the non - endpoint results ($p_{j} > 2$), and \cite{KT} for the endpoint case. See \cite{Tao} for a nice overview of this work. $\Box$\vspace{5mm}

\noindent We will also utilize the maximal Strichartz estimate of \cite{KV2}, which was introduced to provide a new proof of global well - posedness and scattering for the defocusing, three - dimensional energy - critical problem.

\begin{theorem}[Maximal Strichartz estimate]\label{maximalStrichartz}
Suppose that $t, t_{0} \in I$, and

\begin{equation}\label{2.17}
v(t) = \int_{t_{0}}^{t} e^{i(t - \tau) \Delta} F(\tau) d\tau.
\end{equation}

\noindent Then for any $d \geq 3$, $q > \frac{2d}{d - 2}$,

\begin{equation}\label{2.18}
\| \sup_{j} 2^{j(\frac{d}{q} - (d - 2))} \| P_{j} v(t) \|_{L_{x}^{q}(\mathbf{R}^{d})} \|_{L_{t}^{2}(I)} \lesssim \| F \|_{L_{t}^{2} L_{x}^{1}(I \times \mathbf{R}^{d})}.
\end{equation}
\end{theorem}

\noindent \emph{Proof:} This is proved by combining the dispersive estimate $(\ref{2.14})$ with the Sobolev embedding theorem (theorem $\ref{t2.6}$). If $q > \frac{2d}{d - 2}$ then $d(\frac{1}{2} - \frac{1}{q}) > 1$, so

\begin{equation}\label{2.19}
\aligned
2^{j(\frac{d}{q} - (d - 2))} \int_{|t - \tau| > 2^{-2j}} \frac{1}{(t - \tau)^{d(\frac{1}{2} - \frac{1}{q})}} \| P_{j} F(u(\tau)) \|_{L_{x}^{q'}(\mathbf{R}^{d})} d\tau \\ \lesssim \sum_{k \geq 0} 2^{-kd(\frac{1}{2} - \frac{1}{q})} 2^{2j} \int_{|t - \tau| \sim 2^{k} 2^{-2j}} \| F(\tau) \|_{L_{x}^{1}(\mathbf{R}^{d})} d\tau \lesssim_{q} \mathcal M(\| F(\tau) \|_{L_{x}^{1}(\mathbf{R}^{d})})(t).
\endaligned
\end{equation}

\noindent Also by Sobolev embedding

\begin{equation}\label{2.20}
\aligned
2^{j(\frac{d}{q} - (d - 2))} \| \int_{|t - \tau| \leq 2^{-2j}} P_{j} e^{i(t - \tau) \Delta} F(u(\tau)) d\tau \|_{L_{x}^{q}(\mathbf{R}^{d})} \\ \lesssim 2^{2j} \int_{|t - \tau| \leq 2^{-2j}} \| F(\tau) \|_{L_{x}^{1}(\mathbf{R}^{d})} d\tau \lesssim \mathcal M(\| F(\tau) \|_{L_{x}^{1}(\mathbf{R}^{d})})(t).
\endaligned
\end{equation}

\noindent Therefore,

\begin{equation}\label{2.21}
2^{j(\frac{d}{q} - (d - 2))} \| P_{j} v(t) \|_{L_{x}^{q}(\mathbf{R}^{d})} \lesssim \mathcal M(\| F(\tau) \|_{L_{x}^{1}(\mathbf{R}^{d})})(t),
\end{equation}

\noindent so by theorem $\ref{t2.5}$ the proof is complete. $\Box$\vspace{5mm}

\noindent We conclude the section by discussing the double Duhamel trick. This technique was introduced in \cite{CKSTT4} to study the defocusing, energy - critical Schr{\"o}dinger initial value problem when $d = 3$, and in \cite{KV1} for the focusing energy - critical problem for dimensions $d \geq 5$. See also \cite{TaoT1}. The double Duhamel trick is also used to study wave (\cite{Bulut1}, \cite{Bulut2}, \cite{Bulut3}, \cite{KV4}, \cite{KV3}) and KdV (\cite{D9}) problems.\vspace{5mm}

\noindent Suppose $I = [t_{-}, t_{+}]$ and $u$ solves the equation

\begin{equation}\label{2.22}
(i \partial_{t} + \Delta) u = F(t) + G(t).
\end{equation}

\noindent Then by $(\ref{2.11})$, for any $t \in I$,

\begin{equation}\label{2.23}
\aligned
e^{i(t - t_{-}) \Delta} u(t_{-}) - i \int_{t_{-}}^{t} e^{i(t - s_{-}) \Delta} F(s_{-}) ds_{-} - i \int_{t_{-}}^{t} e^{i(t - s_{-}) \Delta} G(s_{-}) ds_{-} = u(t) \\ = e^{i(t - t_{+}) \Delta} u(t_{+}) - i \int_{t_{+}}^{t} e^{i(t - s_{+}) \Delta} F(s_{+}) ds_{+} - i \int_{t_{+}}^{t} e^{i(t - s_{+}) \Delta} G(s_{+}) ds_{+}.
\endaligned
\end{equation}

\noindent Then if $X$ is some Hilbert space, such as $L^{2}(\mathbf{R}^{d})$ or the weighted $L^{2}$ space that we will use in this paper, then by simple linear algebra, for $A + B = A' + B'$,

\begin{equation}\label{2.24}
\langle A + B, A' + B' \rangle \lesssim |A|^{2} + |A'|^{2} + \langle B, B' \rangle,
\end{equation}

\noindent so then

\begin{equation}\label{2.25}
\aligned
\| u(t) \|_{X}^{2} \lesssim \| e^{i(t - t_{-}) \Delta} u(t_{-}) \|_{X}^{2} + \| e^{i(t - t_{+}) \Delta} u(t_{+}) \|_{X}^{2} + \| \int_{t_{-}}^{t} e^{i(t - s_{-}) \Delta} F(s_{-}) ds_{-} \|_{X}^{2} \\ + \| \int_{t_{+}}^{t} e^{i(t - s_{+}) \Delta} F(s_{+}) ds_{+} \|_{X}^{2} + \langle \int_{t_{-}}^{t} e^{i(t - s_{-}) \Delta} G(s_{-}) ds_{-} , e^{i(t - s_{+}) \Delta} G(s_{+}) ds_{+} \rangle_{X}.
\endaligned
\end{equation}

\section{Concentration compactness}
In this section we briefly discuss some concentration compactness results that will be used in this paper.\vspace{5mm}

\noindent \emph{Sketch of the proof of theorem $\ref{t1.9}$:} The reader should consult \cite{KM1} or \cite{KV1} for a complete treatment of the concentration compactness method. Recall that the hypotheses of theorem $\ref{t1.9}$ imply that there exists $E_{\ast} < \| \nabla W \|_{L^{2}(\mathbf{R}^{4})}$ such that

\begin{equation}\label{6.1}
\| u(t) \|_{L_{t}^{\infty} \dot{H}_{x}^{1}(I \times \mathbf{R}^{4})} \leq E_{\ast},
\end{equation}

\noindent where $I$ is the maximal interval of existence for a solution to $(\ref{equation})$. Now let

\begin{equation}\label{6.2}
C(E) = \sup \{ \| u \|_{L_{t,x}^{\frac{2(d + 2)}{d - 2}}(\mathbf{R} \times \mathbf{R}^{d})} : \| u \|_{L_{t}^{\infty} \dot{H}_{x}^{1}(\mathbf{R} \times \mathbf{R}^{d})} \leq E \}.
\end{equation}

\noindent By the results of \cite{CaWe1}, $C(E) \lesssim E$ for $E$ small. Moreover, by a stability result in $d \geq 5$ (see \cite{KV1}) and a simple calculation in dimensions $d = 3, 4$, $C(E)$ is a continuous function of $E$. Therefore, if there exists a non-scattering solution to $(\ref{equation})$ satisfying $(\ref{6.1})$, then by the continuity of $C(E)$, there exists $E_{\ast} < \| \nabla W \|_{L^{2}}$ such that $C(E_{\ast}) = \infty$ and $C(E) < \infty$ for all $E < E_{\ast}$.\vspace{5mm}

\noindent Now take a sequence $u_{n}(t)$ of solutions to $(\ref{equation})$ such that

\begin{equation}\label{6.3}
\| u_{n}(t) \|_{L_{t}^{\infty} \dot{H}_{x}^{1}(\mathbf{R} \times \mathbf{R}^{d})} \nearrow E_{\ast},
\end{equation}

\noindent and

\begin{equation}\label{6.4}
S_{[0, \infty)} (u_{n}) = S_{(-\infty, 0]}(u_{n}) = n.
\end{equation}

\noindent Then \cite{Keraani1} proved that $u_{n}(0)$ can be decomposed into asymptotically decoupling profiles, such that for any $J$,

\begin{equation}\label{6.5}
u_{n}(0) = \sum_{j = 1}^{J} g_{n}^{j} e^{i t_{n}^{j} \Delta} \phi^{j} + w_{n}^{J},
\end{equation}

\noindent where $g_{n}^{j}$ is an element of a group generated by scaling and translation symmetries, $w_{n}^{J}$ represents an error, and the group elements $g_{n}^{j}$ asymptotically decouple. The asymptotic decoupling implies that if $u^{j}(t)$ is the solution to $(\ref{equation})$ with initial data given by $\phi^{j}$, then $(\ref{6.4})$ implies that for one $j_{0}$, $t_{n}^{j_{0}} \rightarrow 0$ and

\begin{equation}\label{6.6}
\| u^{j_{0}}(t) \|_{L_{t}^{\infty} \dot{H}_{x}^{1}(I \times \mathbf{R}^{d})} = E_{\ast},
\end{equation}

\noindent where $I$ is the maximal interval of existence for $u^{j_{0}}(t)$, all other $\phi^{j} = 0$, and

\begin{equation}\label{6.7}
\| u^{j_{0}}(t) \|_{L_{t,x}^{\frac{2(d + 2)}{d - 2}}(I \times \mathbf{R}^{d})} = \infty,
\end{equation}

\noindent where $I$ is the maximal interval of existence of $u^{j_{0}}$. Making the above argument again for a sequence $u^{j_{0}}(t_{n})$, $t_{n} \in I$ shows that $u^{j_{0}}(t_{n})$ has a subsequence that converges in $\dot{H}^{1} / G$, where $G$ is the group of symmetries $g_{n}^{j}$. This proves theorem $\ref{t1.9}$. $\Box$\vspace{5mm}

\noindent Notice that by the Arzela - Ascoli theorem, if $u$ is an almost periodic solution to $(\ref{equation})$, then there exists $x(t) : I \rightarrow \mathbf{R}^{d}$ and $N(t) : I \rightarrow (0, \infty)$, such that for any $\eta > 0$ there exists $C(\eta) < \infty$ such that

\begin{equation}\label{6.8}
\int_{|x - x(t)| > \frac{C(\eta)}{N(t)}} |\nabla u(t,x)|^{2} dx + \int_{|\xi| > C(\eta) N(t)} |\xi|^{2} |\hat{u}(t,\xi)|^{2} d\xi + \int_{|\xi| < \frac{1}{C(\eta)} N(t)} |\xi|^{2} |\hat{u}(t,\xi)|^{2} d\xi < \eta.
\end{equation}

\noindent Moreover, (see \cite{KilVis} for a proof)

\begin{equation}\label{6.9}
|N'(t)| \lesssim N(t)^{3},
\end{equation}

\noindent and

\begin{equation}\label{6.10}
\int_{I} N(t)^{2} dt \lesssim \int_{I} \int |u(t,x)|^{\frac{2(d + 2)}{d - 2}} dx dt \lesssim \int_{I} N(t)^{2} dt + 1.
\end{equation}

\noindent Making use of $(\ref{6.9})$, \cite{KV1} proved theorem $\ref{t1.11}$.\vspace{5mm}

\noindent \emph{Sketch of proof of theorem $\ref{t1.11}$:} Suppose $u(t)$ is an almost periodic solution to $(\ref{equation})$. Then \cite{KV1} showed that one can take a limit of $u(t_{n})$ in $\dot{H}^{1} / G$ and obtain a solution to $(\ref{equation})$ satisfying either

\begin{equation}\label{6.11}
N(t) \geq 1, \hspace{5mm} t \in \mathbf{R}, \hspace{5mm} N(0) = 1,
\end{equation}

\noindent or than $u$ blows up in finite time. However, finite time blowup fails to occur due to concentration compactness and conservation of mass $(\ref{1.3})$. Indeed, by $(\ref{6.10})$, if $u$ blows up in finite time, say at $T = 0$, $N(t) \nearrow \infty$ as $t \searrow 0$. Let $\psi \in C_{0}^{\infty}(\mathbf{R}^{d})$ be a radial function, $\psi = 1$ on $|x| \leq 1$, $\psi$ supported on $|x| \leq 2$. By $(\ref{6.8})$ and H{\"o}lder's inequality, for any $R > 0$,

\begin{equation}\label{6.12}
\lim_{t \searrow 0} \int \psi(\frac{x}{R})^{2} |u(t,x)|^{2} dx = \lim_{t \searrow 0} M_{R}(t) = 0.
\end{equation}

\noindent Moreover, integrating by parts,

\begin{equation}\label{6.13}
\frac{d}{dt} M_{R}(t) \leq \frac{1}{R} \psi'(\frac{x}{R}) \psi(\frac{x}{R}) |\nabla u(t,x)| |u(t,x)| dx \leq \frac{1}{R} M_{R}(t)^{1/2} \| \nabla u(t) \|_{L_{x}^{2}(\mathbf{R}^{d})}.
\end{equation}

\noindent Therefore, $(\ref{6.12})$ combined with the fundamental theorem of calculus and $(\ref{6.13})$ implies that $\int |u(t,x)|^{2} dx = 0$ for any $t > 0$. However, this implies $u \equiv 0$, which contradicts $u$ blowing up in finite time. $\Box$

\begin{theorem}\label{t6.1}
If $u(t)$ is an almost periodic solution to $(\ref{equation})$ satisfying $N(t) \geq 1$, then

\begin{equation}\label{6.14}
\| u(t) \|_{L_{t}^{\infty} L_{x}^{3}(\mathbf{R} \times \mathbf{R}^{4})} < \infty.
\end{equation}
\end{theorem}

\noindent \textbf{Remark:} This is an endpoint of a more general result of \cite{MMZ}.\vspace{5mm}

\noindent \textbf{Remark:} From now on since we are considering an almost periodic solution $u$ to $(\ref{equation})$, let $A \lesssim B$ denote $A \leq C(u) B$.\vspace{5mm}

\noindent \emph{Proof:} By the Duhamel formula $(\ref{2.11})$, for any $t_{0} \in \mathbf{R}$,

\begin{equation}\label{6.15}
u(t) = e^{i(t - t_{0}) \Delta} u(t_{0}) - i \int_{t_{0}}^{t} e^{i(t - \tau) \Delta} F(u(\tau)) d\tau.
\end{equation}

\noindent Now by $(\ref{6.8})$, for a fixed $t$,

\begin{equation}\label{6.16}
e^{i(t - t_{0}) \Delta} u(t_{0}) \rightharpoonup 0
\end{equation}

\noindent weakly in $\dot{H}^{1}$ as $t_{0} \rightarrow \pm \infty$. Indeed, this follows from $(\ref{6.8})$ if $N(t_{0}) \nearrow +\infty$ or $\searrow 0$. If $N(t_{0}) \sim 1$ then this follows by combining $(\ref{6.8})$ and the dispersive estimate $(\ref{2.14})$. Therefore, for any $j \in \mathbf{Z}$,

\begin{equation}\label{6.17}
\| P_{j} u(t) \|_{L_{x}^{\infty}(\mathbf{R}^{4})} \lesssim \lim_{t_{0} \rightarrow \infty} \| \int_{t_{0}}^{t} e^{i(t - \tau) \Delta} P_{j} F(u(\tau)) d\tau \|_{L_{x}^{\infty}(\mathbf{R}^{4})}.
\end{equation}

\noindent Now by Sobolev embedding $(\ref{2.7})$ and $(\ref{2.12})$,

\begin{equation}\label{6.18}
\| \int_{t - 2^{-2j}}^{t} e^{i(t - \tau) \Delta} P_{j} F(u(\tau)) d\tau \|_{L_{x}^{\infty}(\mathbf{R}^{4})} \lesssim 2^{2j} \| u \|_{L_{t}^{\infty} L_{x}^{3}(\mathbf{R} \times \mathbf{R}^{4})}^{3}.
\end{equation}

\noindent Also by the dispersive estimate $(\ref{2.14})$,

\begin{equation}\label{6.19}
\aligned
\| \int_{t_{0}}^{t - 2^{-2j}} e^{i(t - \tau) \Delta} P_{j} F(u(\tau)) d\tau \|_{L_{x}^{\infty}(\mathbf{R}^{4})} \\ \lesssim \| u \|_{L_{t}^{\infty} L_{x}^{3}(\mathbf{R} \times \mathbf{R}^{4})}^{3} \int_{t > 2^{-2j}} \frac{1}{t^{2}} dt \lesssim 2^{2j} \| u \|_{L_{t}^{\infty} L_{x}^{3}(\mathbf{R} \times \mathbf{R}^{4})}^{3}.
\endaligned
\end{equation}

\noindent Then by $(\ref{6.8})$, if $j_{0}(\eta)$ is the largest integer such that $2^{j_{0}} \leq \frac{1}{C(\eta)}$, then by $N(t) \geq 1$, $(\ref{6.18})$ and $(\ref{6.19})$

\begin{equation}\label{6.20}
\aligned
\| P_{\leq j_{0}(\eta)} u(t) \|_{L_{x}^{3}(\mathbf{R}^{4})}^{3} \lesssim \sum_{k_{1} \leq k_{2} \leq k_{3} \leq j_{0}(\eta)} \| P_{k_{1}} u(t) \|_{L_{x}^{\infty}(\mathbf{R}^{4})} \| P_{k_{2}} u(t) \|_{L_{x}^{2}(\mathbf{R}^{4})} \| P_{k_{3}} u(t) \|_{L_{x}^{2}(\mathbf{R}^{4})} \\
\lesssim \| u \|_{L_{t}^{\infty} L_{x}^{3}(\mathbf{R} \times \mathbf{R}^{4})}^{3} \sum_{k_{1} \leq k_{2} \leq k_{3} \leq j_{0}(\eta)} 2^{2k_{1}} 2^{-k_{2}} 2^{-k_{3}} \| \nabla P_{k_{2}} u(t) \|_{L_{x}^{2}(\mathbf{R}^{4})} \| \nabla P_{k_{3}} u(t) \|_{L_{x}^{2}(\mathbf{R}^{4})} \\ \lesssim \eta^{2} \| u \|_{L_{t}^{\infty} L_{x}^{3}(\mathbf{R} \times \mathbf{R}^{4})}^{3}.
\endaligned
\end{equation}

\noindent Meanwhile, by Bernstein's inequality, since $\dot{H}^{2/3}(\mathbf{R}^{4}) \subset L_{x}^{3}(\mathbf{R}^{4})$,

\begin{equation}\label{6.21}
\| P_{> j_{0}} u(t) \|_{L_{x}^{3}(\mathbf{R}^{4})}^{3} \lesssim 2^{-j_{0}},
\end{equation}

\noindent so

\begin{equation}\label{6.22}
\| u(t) \|_{L_{t}^{\infty} L_{x}^{3}(\mathbf{R} \times \mathbf{R}^{4})} \lesssim 2^{-j_{0}(\eta)/3} + \eta^{2/3} \| u(t) \|_{L_{t}^{\infty} L_{x}^{3}(\mathbf{R} \times \mathbf{R}^{4})}.
\end{equation}

\noindent $\Box$\vspace{5mm}

\section{Long Time Strichartz estimate}

\noindent Now we prove a long time Strichartz estimate. Long time Strichartz estimates were introduced in \cite{D2} to study the mass - criticial nonlinear Schr{\"o}dinger equation in dimensions $d \geq 3$. Subsequently, \cite{Visan1} and \cite{KV2} utilized long - time Strichartz estimates for the defocusing, energy - critical nonlinear Schr{\"o}dinger problem in dimensions $d = 4$ and $d = 3$ respectively. The long - time Strichartz estimates of $d = 3$ relied on the maximal Strichartz estimates. Here we prove the long - time Strichartz estimates of \cite{KV2}, modified to the case when $d = 4$. In this case the crucial quantity is

\begin{equation}\label{3.2}
K = \int_{I} N(t)^{-2} dt.
\end{equation}

\noindent The quantity $\int_{I} \frac{1}{N(t)} dt$ was quite useful in the defocusing case since it scaled like the interaction Morawetz estimates of \cite{CKSTT2}, \cite{CKSTT4}, and \cite{TVZ}. However, in the focusing case there is no such estimate. On the other hand, $(\ref{1.19})$ is a lower bound for

\begin{equation}\label{1.20}
\int_{I} \int |u(t,x)|^{2} dx dt,
\end{equation}

\noindent which allows us to prove some very useful interaction Morawetz estimates.\vspace{5mm}

\begin{theorem}[Long time Strichartz estimate]\label{t3.2}
For any $j$,

\begin{equation}\label{3.14}
(\sum_{k \leq j} \| \nabla u_{j} \|_{L_{t}^{2} L_{x}^{4}(I \times \mathbf{R}^{4})}^{2})^{1/2} + 2^{2j} \| \sup_{k \geq j} 2^{-2k} \| u_{k}(t) \|_{L_{x}^{\infty}(\mathbf{R}^{4})} \|_{L_{t}^{2}(I)} \lesssim (1 + K 2^{4j})^{1/2}.
\end{equation}
\end{theorem}

\noindent \emph{Proof:} By Sobolev embedding and Bernstein's inequality

\begin{equation}\label{3.15}
2^{4j} (\sum_{k \geq j} 2^{-4k} \| e^{i(t - t_{0}) \Delta} P_{k} u(t_{0}) \|_{L_{t}^{2} L_{x}^{\infty}(I \times \mathbf{R}^{4})}^{2}) \lesssim 2^{2j} \| P_{> j} u(t_{0}) \|_{L_{x}^{2}(\mathbf{R}^{4})}^{2} \lesssim \| u(t_{0}) \|_{\dot{H}^{1}(\mathbf{R}^{4})}^{2} \lesssim 1,
\end{equation}

\noindent and

\begin{equation}\label{3.15.1}
\sum_{k \leq j} \| \nabla e^{i(t - t_{0}) \Delta} P_{k} u(t_{0}) \|_{L_{t}^{2} L_{x}^{4}(I \times \mathbf{R}^{4})}^{2} \lesssim \| \nabla u(t_{0}) \|_{L_{x}^{2}(\mathbf{R}^{4})}^{2} \lesssim 1.
\end{equation}

\noindent Now let

\begin{equation}\label{3.16}
\aligned
\| u \|_{Y(I \times \mathbf{R}^{4})} = \sup_{j} 2^{2j} (1 + K 2^{4j})^{-1/2} \| \sup_{k \geq j} 2^{-2k} \| u_{j}(t) \|_{L_{x}^{\infty}(\mathbf{R}^{4})} \|_{L_{t}^{2}(I)} \\ + \sup_{j} (1 + K 2^{4j})^{-1/2} (\sum_{k \leq j} 2^{2k} \| u_{k}(t) \|_{L_{t}^{2} L_{x}^{4}(I \times \mathbf{R}^{4})}^{2})^{1/2}.
\endaligned
\end{equation}

\noindent By conservation of energy and Bernstein's inequality,

\begin{equation}\label{3.17}
\aligned
\| P_{\leq cN(t)} P_{\geq j} u \|_{L_{t}^{6} L_{x}^{3}(I \times \mathbf{R}^{4})}^{3} \\ \lesssim \| \sum_{j \leq k_{1} \leq k_{2} \leq k_{3}} \| P_{\leq cN(t)} u_{k_{1}} \|_{L_{x}^{\infty}(\mathbf{R}^{4})} \| P_{\leq cN(t)} u_{k_{2}} \|_{L_{x}^{2}(\mathbf{R}^{4})} \| P_{\leq cN(t)} u_{k_{3}} \|_{L_{x}^{2}(\mathbf{R}^{4})} \|_{L_{t}^{2}(I)} \\ \lesssim \eta^{2} \| \sup_{k \geq j} 2^{-2k} \| u_{k}(t) \|_{L_{x}^{\infty}(\mathbf{R}^{4})} \|_{L_{t}^{2}(I)}.
\endaligned
\end{equation}

\noindent Also, by Bernstein's inequality

\begin{equation}\label{3.18}
\| P_{\geq cN(t)} u \|_{L_{t}^{6} L_{x}^{3}(I \times \mathbf{R}^{4})}^{3} \lesssim (\int_{I} c^{-2} N(t)^{-2} dt)^{1/2} \lesssim c^{-1} K^{1/2}.
\end{equation}

\noindent Therefore,

\begin{equation}\label{3.19}
\| u_{\geq j}^{3} \|_{L_{t}^{2} L_{x}^{1}(I \times \mathbf{R}^{4})} \lesssim c^{-1} K^{1/2} + \eta^{2} (1 + 2^{4j} K)^{1/2} 2^{-2j} \| u \|_{Y(I \times \mathbf{R}^{4})},
\end{equation}

\noindent and by theorem $\ref{maximalStrichartz}$ and Sobolev embedding,

\begin{equation}\label{3.20}
\aligned
\| \sup_{k \geq j} \| P_{k} \int_{t_{0}}^{t} e^{i(t - \tau) \Delta} F(u_{\geq j}) d\tau \|_{L_{x}^{\infty}(\mathbf{R}^{4})} \|_{L_{t}^{2}(I)} \\ \lesssim c^{-1} K^{1/2} + \eta^{2} (1 + 2^{4j} K)^{1/2} 2^{-2j} \| u \|_{Y(I \times \mathbf{R}^{4})},
\endaligned
\end{equation}

\noindent and

\begin{equation}\label{3.21}
\aligned
(\sum_{k \leq j} 2^{2k} \| P_{k} \int_{t_{0}}^{t} e^{i(t - \tau) \Delta} F(u_{\geq j}) d\tau \|_{L_{t}^{2} L_{x}^{4}(I \times \mathbf{R}^{4})}^{2})^{1/2} \\ \lesssim c^{-1} 2^{2j} K^{1/2} + \eta^{2} (1 + 2^{4j} K)^{1/2} \| u \|_{Y(I \times \mathbf{R}^{4})}.
\endaligned
\end{equation}

\noindent Next, by the Littlewood - Paley theorem and Sobolev embedding,

\begin{equation}\label{3.22}
\aligned
\| u_{\leq j} \|_{L_{t,x}^{6}(I \times \mathbf{R}^{4})} \lesssim \| \nabla u_{\leq j} \|_{L_{t}^{6} L_{x}^{12/5}(I \times \mathbf{R}^{4})} \\ \lesssim (\sum_{k \leq j} 2^{2k} \| u_{k} \|_{L_{t}^{2} L_{x}^{4}(I \times \mathbf{R}^{4})}^{2})^{1/6} \| \nabla u \|_{L_{t}^{\infty} L_{x}^{2}(I \times \mathbf{R}^{4})}^{2/3} \lesssim (1 + 2^{4j} K)^{1/6} \| u \|_{Y(I \times \mathbf{R}^{4})}^{1/3}.
\endaligned
\end{equation}

\noindent Therefore by Sobolev embedding, $(\ref{3.19})$ and $(\ref{3.22})$,

\begin{equation}\label{3.23}
\aligned
\| (u_{\geq j}^{2} u_{\leq j}) \|_{L_{t}^{2} L_{x}^{4/3}(I \times \mathbf{R}^{4})} \lesssim \| u_{\geq j} \|_{L_{t}^{6} L_{x}^{3}(I \times \mathbf{R}^{4})}^{2} \| u_{\leq j} \|_{L_{t}^{6} L_{x}^{12}(I \times \mathbf{R}^{4})} \\ \lesssim 2^{j/3} (c^{-1} K^{1/2} + \eta^{2} 2^{-2j} (1 + 2^{4j} K)^{1/2} \| u \|_{Y(I \times \mathbf{R}^{4})})^{2/3} (1 + 2^{4j} K)^{1/6} \| u \|_{Y(I \times \mathbf{R}^{4})}^{1/3} \\ \lesssim 2^{j/3} c^{-2/3} K^{1/3} (1 + 2^{4j} K)^{1/6} \| u \|_{Y(I \times \mathbf{R}^{4})}^{1/3} + 2^{-j} \eta^{4/3} (1 + 2^{4j} K)^{1/2} \| u \|_{Y(I \times \mathbf{R}^{4})}.
\endaligned
\end{equation}

\noindent $(\ref{3.23})$ implies, by Sobolev embedding, that

\begin{equation}\label{3.24}
\aligned
\| \sup_{k \geq j} 2^{-2k} \| P_{k} \int_{t_{0}}^{t} e^{i(t - \tau) \Delta} O(u_{\geq j}^{2} u_{\leq j}) d\tau \|_{L_{x}^{\infty}(\mathbf{R}^{4})} \|_{L_{t}^{2}(I)} \\ \lesssim 2^{-j} \| \int_{t_{0}}^{t} e^{i(t - \tau) \Delta} O(u_{\geq j}^{2} u_{\leq j}) d\tau \|_{L_{t}^{2} L_{x}^{4}(I \times \mathbf{R}^{4})} \\ \lesssim 2^{-2j/3} c^{-2/3} K^{1/3} (1 + 2^{4j} K)^{1/6} \| u \|_{Y(I \times \mathbf{R}^{4})}^{1/3} \\ + 2^{-2j} \eta^{4/3} (1 + 2^{4j} K)^{1/2} \| u \|_{Y(I \times \mathbf{R}^{4})},
\endaligned
\end{equation}

\noindent and by Strichartz estimates

\begin{equation}\label{3.25}
\aligned
(\sum_{k \leq j} 2^{2k} \| P_{k} \int_{t_{0}}^{t} e^{i(t - \tau) \Delta} O(u_{\geq j}^{2} u_{\leq j}) d\tau \|_{L_{t}^{2} L_{x}^{4}(I \times \mathbf{R}^{4})}^{2})^{1/2} \lesssim 2^{j} \| u_{\geq j}^{2} u_{\leq j} \|_{L_{t}^{2} L_{x}^{4/3}(I \times \mathbf{R}^{4})} \\ \lesssim c^{-2/3} 2^{4j/3} K^{1/3} (1 + 2^{4j} K)^{1/6} \| u \|_{Y(I \times \mathbf{R}^{4})}^{1/3} + \eta^{4/3} (1 + 2^{4j} K)^{1/2} \| u \|_{Y(I \times \mathbf{R}^{4})}.
\endaligned
\end{equation}

\noindent Next, by Sobolev embedding,

\begin{equation}\label{3.26}
\aligned
\| (P_{\leq c N(t)} u_{\leq j})^{2} \|_{L_{t}^{2} L_{x}^{4}(I \times \mathbf{R}^{4})} \lesssim \| \nabla u_{\leq j} \|_{L_{t}^{2} L_{x}^{4}(I \times \mathbf{R}^{4})} \| u_{\leq cN(t)} \|_{L_{t}^{\infty} L_{x}^{4}(I \times \mathbf{R}^{4})} \\ \lesssim \eta (1 + 2^{4j} K)^{1/2} \| u \|_{Y(I \times \mathbf{R}^{4})},
\endaligned
\end{equation}

\noindent and by Bernstein's inequality and Sobolev embedding

\begin{equation}\label{3.27}
\aligned
\| (P_{\leq c N(t)} u_{\leq j})^{2} \|_{L_{t}^{2} L_{x}^{4}(I \times \mathbf{R}^{4})} \lesssim 2^{2j} (\int \| u_{> cN(t)} \|_{L_{x}^{2}(\mathbf{R}^{4})}^{2} dt)^{1/2} \lesssim c^{-1} K^{1/2} 2^{2j}.
\endaligned
\end{equation}

\noindent $(\ref{3.26})$ and $(\ref{3.27})$ imply that

\begin{equation}\label{3.28}
\aligned
\| \nabla (u_{\leq j}^{2} u_{> j}) \|_{L_{t}^{2} L_{x}^{4/3}(I \times \mathbf{R}^{4})} + \| \nabla u_{\leq j}^{3} \|_{L_{t}^{2} L_{x}^{4/3}(I \times \mathbf{R}^{4})} \\ \lesssim \| \nabla u \|_{L_{t}^{\infty} L_{x}^{2}(I \times \mathbf{R}^{4})} \| u_{\leq j}^{2} \|_{L_{t}^{2} L_{x}^{4}(I \times \mathbf{R}^{4})} \lesssim c^{-1} K^{1/2} 2^{2j} + \eta (1 + 2^{4j} K)^{1/2} \| u \|_{Y(I \times \mathbf{R}^{4})}. 
\endaligned
\end{equation}

\noindent Therefore,

\begin{equation}\label{3.29}
\aligned
(\sum_{k \leq j} 2^{2k} \| P_{k} \int_{t_{0}}^{t} e^{i(t - \tau) \Delta} O(u_{\leq j}^{2} u) d\tau \|_{L_{t}^{2} L_{x}^{4}(I \times \mathbf{R}^{4})})^{1/2} \\ \lesssim c^{-1} K^{1/2} 2^{2j} + \eta (1 + 2^{4j} K)^{1/2} \| u \|_{Y(I \times \mathbf{R}^{4})}.
\endaligned
\end{equation}

\noindent Also, by Sobolev embedding, Bernstein's inequality, and $(\ref{3.29})$,

\begin{equation}\label{3.30}
\aligned
\| \sup_{k \geq j} 2^{-2k} \| P_{k} \int_{t_{0}}^{t} e^{i(t - \tau) \Delta} O(u_{\leq j}^{2} u) d\tau \|_{L_{x}^{\infty}(\mathbf{R}^{4})} \|_{L_{t}^{2}(I)} \\ \lesssim c^{-1} K^{1/2} + \eta 2^{-2j} (1 + 2^{4j} K)^{1/2} \| u \|_{Y(I \times \mathbf{R}^{4})}.
\endaligned
\end{equation}

\noindent Therefore, combining $(\ref{3.20})$, $(\ref{3.21})$, $(\ref{3.24})$, $(\ref{3.25})$, $(\ref{3.29})$, and $(\ref{3.30})$,

\begin{equation}\label{3.31}
\| u \|_{Y(I \times \mathbf{R}^{4})} \lesssim c(\eta)^{-1} + \eta \| u \|_{Y(I \times \mathbf{R}^{4})}.
\end{equation}

\noindent Choosing $\eta > 0$ sufficiently small and then $c(\eta) > 0$ sufficiently small, the proof of theorem $\ref{t3.2}$ is complete. $\Box$\vspace{5mm}

\noindent Now, armed with the long time Strichartz estimate we can rule out the rapid frequency cascade scenario.

\begin{theorem}\label{t3.3}
If $u$ is an almost periodic solution to $(\ref{equation})$ on $\mathbf{R}$, $\int N(t)^{-2} dt = K < \infty$, then $u \equiv 0$.
\end{theorem}

\noindent \emph{Proof:} Let $k_{0}$ be the integer closest to $k$, where $2^{k} = K^{-1/4}$. Choose $j \leq k_{0}$. Then

\begin{equation}\label{3.32}
\| \nabla P_{\leq j} u(t) \|_{L_{x}^{2}(\mathbf{R}^{4})} \lesssim \| \nabla P_{\leq j} u(-T) \|_{L_{x}^{2}(\mathbf{R}^{4})} + \| \nabla P_{\leq j} \int_{-T}^{t} e^{i(t - \tau) \Delta} F(u(\tau)) d\tau \|_{L_{x}^{2}(\mathbf{R}^{4})}.
\end{equation}

\noindent Also choose $j < j_{0}(\eta)$ so that

\begin{equation}\label{3.33}
\| P_{\leq j} u(t) \|_{L_{t}^{\infty} \dot{H}_{x}^{1}(\mathbf{R} \times \mathbf{R}^{4})} \leq \eta.
\end{equation}

\noindent Then

\begin{equation}\label{3.33}
\| \nabla F(u_{\leq j}) \|_{L_{t}^{2} L_{x}^{4/3}([-T, T] \times \mathbf{R}^{4})} \lesssim \eta^{2} \| \nabla u_{\leq j} \|_{L_{t}^{2} L_{x}^{4}([-T, T] \times \mathbf{R}^{4})}.
\end{equation}

\noindent Next, by Sobolev embedding,

\begin{equation}\label{3.34}
\aligned
\| \nabla P_{\leq j} O(u_{\leq j}^{2} u_{> j} \|_{L_{t}^{2} L_{x}^{4/3}([-T, T] \times \mathbf{R}^{4})} \\ \lesssim 2^{j} \| u_{j \leq \cdot \leq k_{0}} \|_{L_{t}^{\infty} L_{x}^{2}(\mathbf{R} \times \mathbf{R}^{4})} \| \nabla u_{\leq j} \|_{L_{t}^{2} L_{x}^{4}([-T, T] \times \mathbf{R}^{4})} \| \nabla u_{\leq j} \|_{L_{t}^{\infty} L_{x}^{2}(\mathbf{R} \times \mathbf{R}^{4})}.
\endaligned
\end{equation}

\noindent Finally, by Bernstein's inequality and Sobolev embedding,

\begin{equation}\label{3.35}
\aligned
\| \nabla P_{\leq j} O(u_{j < \cdot < k_{0}}^{2} u) \|_{L_{t}^{2} L_{x}^{4/3}([-T, T] \times \mathbf{R}^{4})} \\ \lesssim 2^{2j} \sum_{j \leq k_{1} \leq k_{2} \leq k_{0}} \| P_{k_{1}} u \|_{L_{t}^{2} L_{x}^{4}([-T, T] \times \mathbf{R}^{4})} \| P_{k_{2}} u \|_{L_{t}^{\infty} L_{x}^{2}([-T, T] \times \mathbf{R}^{4})} \| u \|_{L_{t}^{\infty} L_{x}^{4}([-T, T] \times \mathbf{R}^{4})},
\endaligned
\end{equation}

\noindent and by theorem $\ref{t3.2}$ and Sobolev embedding,

\begin{equation}\label{3.36}
\| \nabla P_{\leq j} O(u_{\geq k_{0}}^{3}) \|_{L_{t}^{2} L_{x}^{4/3}([-T, T] \times \mathbf{R}^{4})} \lesssim 2^{2j} \| u_{\geq k_{0}} \|_{L_{t}^{6} L_{x}^{3}([-T, T] \times \mathbf{R}^{4})}^{3} \lesssim 2^{2j} K^{1/2},
\end{equation}

\noindent and

\begin{equation}\label{3.37}
\aligned
\| \nabla P_{\leq j} O(u_{\geq k_{0}}^{2} u_{\leq k_{0}}) \|_{L_{t}^{2} L_{x}^{4/3}([-T, T] \times \mathbf{R}^{4})} \\ \lesssim 2^{3j/2} \| u_{\leq k_{0}} \|_{L_{t,x}^{6} ([-T, T] \times \mathbf{R}^{4})} \| u_{> k_{0}} \|_{L_{t}^{6} L_{x}^{3}([-T, T] \times \mathbf{R}^{4})} \lesssim 2^{3j/2} K^{1/2}.
\endaligned
\end{equation}

 \noindent Therefore, by $(\ref{3.32})$ - $(\ref{3.37})$, by $(\ref{6.8})$ and $(\ref{3.33})$, for any $T$,

\begin{equation}\label{3.38}
\aligned
\| \nabla P_{\leq j} u(t) \|_{L_{t}^{2} L_{x}^{4}([-T, T] \times \mathbf{R}^{4})} \\ \lesssim \| \nabla P_{\leq j} u(-T) \|_{L_{x}^{2}(\mathbf{R}^{4})} + \eta^{2} \| \nabla P_{\leq j} u \|_{L_{t}^{2} L_{x}^{4}([-T, T] \times \mathbf{R}^{4})} \\ + \sum_{j \leq l \leq k_{0}} 2^{2(j - l)} \| \nabla P_{\leq l} u(t) \|_{L_{t}^{2} L_{x}^{4}([-T, T] \times \mathbf{R}^{4})} + 2^{3j/2} K^{1/2}.
\endaligned
\end{equation}

\noindent By theorem $\ref{t3.2}$, for $l \leq k_{0}$,

\begin{equation}\label{3.39}
\| \nabla P_{l} u \|_{L_{t}^{2} L_{x}^{4}([-T, T] \times \mathbf{R}^{4})} \lesssim 1,
\end{equation}

\noindent uniformly in $T$. Also, since $\int_{\mathbf{R}} N(t)^{-2} dt = K < \infty$, $N(-T) \nearrow +\infty$ as $T \nearrow +\infty$, so

\begin{equation}\label{3.40}
\inf_{T} \| \nabla P_{\leq j} u(-T) \|_{L_{x}^{2}(\mathbf{R}^{4})} = 0.
\end{equation}

\noindent Therefore, by induction, starting with $(\ref{3.39})$ for $j_{0}(\eta) \leq l \leq k_{0}$, $(\ref{3.38})$ implies

\begin{equation}\label{3.41}
\| \nabla P_{\leq j} u(t) \|_{L_{t}^{2} L_{x}^{4}(\mathbf{R} \times \mathbf{R}^{4})} \lesssim K^{1/2} 2^{3j/2}.
\end{equation}

\noindent Also since

\begin{equation}\label{3.42}
\aligned
\| \nabla P_{\leq j} u(t) \|_{L_{t}^{\infty} L_{x}^{2}([-T,T] \times \mathbf{R}^{4})} \lesssim \| \nabla P_{\leq j} u(-T) \|_{L_{x}^{2}(\mathbf{R}^{4})} + \eta^{2} \| \nabla P_{\leq j} u \|_{L_{t}^{2} L_{x}^{4}([-T, T] \times \mathbf{R}^{4})} \\ + \sum_{j \leq l \leq k_{0}} 2^{2(j - l)} \| \nabla P_{\leq l} u(t) \|_{L_{t}^{2} L_{x}^{4}([-T, T] \times \mathbf{R}^{4})} + 2^{3j/2} K^{1/2},
\endaligned
\end{equation}

\begin{equation}\label{3.43}
\| \nabla P_{\leq j} u(t) \|_{L_{t}^{\infty} L_{x}^{2}(\mathbf{R} \times \mathbf{R}^{4})} \lesssim K^{1/2} 2^{3j/2}.
\end{equation}

\noindent In particular this implies 

\begin{equation}\label{3.44}
\| u(t) \|_{H_{x}^{-1/4}(\mathbf{R}^{4})} \lesssim K^{5/12}.
\end{equation}

\noindent Then by Bernstein's inequality, interpolation, $(\ref{6.8})$, and $(\ref{3.44})$,

\begin{equation}\label{3.45}
\aligned
\| u(t) \|_{L_{x}^{2}(\mathbf{R}^{4})} \lesssim \| P_{\leq \frac{1}{C(\eta)} N(t)} u(t) \|_{H_{x}^{-1/4}(\mathbf{R}^{4})}^{4/5} \| P_{\leq \frac{1}{C(\eta)} N(t)} u(t) \|_{\dot{H}_{x}^{1}(\mathbf{R}^{4})}^{1/5} \\ + \| P_{\geq \frac{1}{C(\eta)} N(t)} u(t) \|_{L_{x}^{2}(\mathbf{R}^{4})} \lesssim K^{2/3} \eta^{1/5} + \frac{C(\eta)}{N(t)}.
\endaligned
\end{equation}

\noindent Then $\int_{\mathbf{R}} N(t)^{-2} dt = K$ combined with $(\ref{6.9})$ implies that $N(t) \nearrow +\infty$ as $t \nearrow \infty$, so we can choose $\eta(t) \searrow 0$, possibly very slowly, such that

\begin{equation}\label{3.46}
\| u(t) \|_{L_{x}^{2}(\mathbf{R}^{4})} \rightarrow 0.
\end{equation}

\noindent This implies $u \equiv 0$. $\Box$

\section{Soliton}
Now we turn to the case when $\int_{\mathbf{R}} N(t)^{-2} dt = \infty$. We begin by excluding the soliton, the case when $N(t) = 1$ for all $t \in \mathbf{R}$. To do this we utilize an interaction Morawetz estimate. As in \cite{CKSTT4} and \cite{KV2} the interaction Morawetz estimate utilizes an integral of mass estimate.

\begin{lemma}\label{l4.1}
Suppose that $\psi \in C_{0}^{\infty}(\mathbf{R}^{4})$, $\psi = 1$ for $|x| \leq 1$, and $\psi$ is supported on $|x| \leq 2$. For any $1 \leq R \leq K^{1/5}$, where $\int_{I} N(t)^{-2} dt = \int_{I} 1 dt = K$,

\begin{equation}\label{4.1}
 \int_{I} \int \int |u(t,y)|^{2} \psi(\frac{x - y}{R}) [|\nabla u(t,x)|^{2} + |u(t,x)|^{4}] dx dy dt \lesssim K \ln(R).
\end{equation}
\end{lemma}

\noindent \emph{Proof:} This is proved using the double Duhamel trick. See \cite{MMZ} for a similar result. Both here and in \cite{MMZ} we have a logarithmic - type failure. Suppose $I = [t_{-}, t_{+}]$.\vspace{5mm}

\noindent Let $P_{h} = P_{\geq K^{-1/4}}$, $P_{h} + P_{l} = 1$. By Duhamel's principle,

\begin{equation}\label{4.2}
\aligned
u_{h}(t) = e^{i(t - t_{-}) \Delta} u_{h}(t_{-}) + \int_{t_{-}}^{t} e^{i(t - s_{-}) \Delta} P_{h} O(u_{l} u^{2}) ds_{-} + \int_{t_{-}}^{t} e^{i(t - s_{-}) \Delta} P_{h} F(u_{h}) ds_{-} \\
= e^{i(t - t_{+}) \Delta} u_{h}(t_{+}) + \int_{t_{+}}^{t} e^{i(t - s_{+}) \Delta} P_{h} O(u_{l} u^{2}) ds_{+} + \int_{t_{+}}^{t} e^{i(t - s_{+}) \Delta} P_{h} F(u_{h}) ds_{+}.
\endaligned
\end{equation}

\noindent Now for a fixed $x \in \mathbf{R}^{4}$ define the inner product

\begin{equation}\label{4.3}
\langle f, g \rangle_{x} = \int \psi(\frac{x - y}{R}) f(y) \overline{g(y)} dy.
\end{equation}

\noindent Now use $(\ref{2.23})$, $(\ref{2.24})$, and $(\ref{2.25})$. Let

\begin{equation}\label{4.4}
\aligned
A = e^{i(t - t_{-}) \Delta} u_{h}(t_{-}) - i \int_{t_{-}}^{t} e^{i(t - s_{-}) \Delta} O(u_{l} u^{2}) ds_{-} \\ - i \int_{t_{-}}^{t - R^{2}} e^{i(t - s_{-}) \Delta} F(u_{h}) ds_{-} - i \int_{t - 1}^{t} e^{i(t - s_{-}) \Delta} F(u_{h}) ds_{-}, \\
A' = e^{i(t - t_{+}) \Delta} u_{h}(t_{+}) - i \int_{t_{+}}^{t} e^{i(t - s_{+}) \Delta} O(u_{l} u^{2}) ds_{+} \\ - i \int_{t_{+}}^{t + R^{2}} e^{i(t - s_{+}) \Delta} F(u_{h}) ds_{+} - i \int_{t + 1}^{t} e^{i(t - s_{+}) \Delta} F(u_{h}) ds_{+},
\endaligned
\end{equation}

\noindent and

\begin{equation}\label{4.5}
\aligned
B = -i  \int_{t - R^{2}}^{t - 1} e^{i(t - s_{-}) \Delta} F(u_{h}) ds_{-}, \\
B' = - i \int_{t + R^{2}}^{t + 1} e^{i(t - s_{+}) \Delta} F(u_{h}) ds_{+}.
\endaligned
\end{equation}

\noindent By $(\ref{3.23})$, $(\ref{3.28})$, Strichartz estimates, Bernstein's inequality, and H{\"o}lder's inequality,

\begin{equation}\label{4.6}
\int_{I} \int \psi(\frac{x - y}{R}) [|e^{i(t - t_{-}) \Delta} u(t_{-})(y)|^{2} + |\int_{t_{-}}^{t} e^{i(t - s_{-}) \Delta} O(u_{l} u^{2})(y) ds_{-}|^{2} ] dy dt \lesssim R^{2} K^{1/2}.
\end{equation}

\noindent By $(\ref{2.14})$ and H{\"o}lder's inequality,

\begin{equation}\label{4.7}
\aligned
 \int_{t_{-}}^{t - R^{2}} \int \psi(\frac{x - y}{R}) |e^{i(t - s_{-}) \Delta} F(u_{h})(s_{-})(y)|^{2} dy \\ \lesssim R^{4} (\int_{t - s_{-} > R^{2}} \frac{1}{(t - s_{-})^{2}} \| F(u_{h})(s_{-}) \|_{L_{x}^{1}} ds_{-})^{2}  \\
\lesssim (\sum_{j \geq 0} \frac{1}{2^{j} R^{2}} \int_{|t - s_{-}| \sim 2^{j} R^{2}} \| F(u_{h})(s_{-}) \|_{L_{x}^{1}(\mathbf{R}^{4})} ds_{-})^{2} \lesssim (\mathcal M(\| F(s) \|_{L_{x}^{1}(\mathbf{R}^{4})})(t))^{2}.
\endaligned
\end{equation}

\noindent Finally, by H{\"o}lder's inequality in time and Sobolev embedding,

\begin{equation}\label{4.8}
\| u_{h}^{3} \|_{L_{t}^{2} L_{x}^{4/3}([t - 1, t] \times \mathbf{R}^{2})} \lesssim 1,
\end{equation}

\noindent so by theorem $\ref{t2.8}$ and $(\ref{4.6})$ - $(\ref{4.8})$, 

\begin{equation}\label{4.9}
|A|^{2} \lesssim 1 + \mathcal M(\| u_{h} \|_{L_{x}^{3}(\mathbf{R}^{4})}^{3})(t)^{2} + a(t)^{2} R^{2},
\end{equation}

\noindent where $\int a(t)^{2} dt \lesssim K^{1/2}$. By an identical calculation

\begin{equation}\label{4.10}
|A'|^{2} \lesssim 1 + \mathcal M(\| u_{h} \|_{L_{x}^{3}(\mathbf{R}^{4})}^{3})(t)^{2} + a(t)^{2} R^{2}.
\end{equation}

\noindent To compute $\langle B, B' \rangle_{x}$, we compute the kernel of $e^{i(t - s_{-}) \Delta} \psi(\frac{x - y}{R}) e^{i (s_{+} - t) \Delta}$. Since $t$ is fixed, to simplify notation let $x = 0$, $s = s_{+} - t$ and $t = t - s_{-}$. Then the kernel of $e^{it \Delta} \psi(\frac{y}{R}) e^{is \Delta}$ is given by

\begin{equation}\label{4.11}
K(s,t; y, z) = \frac{C}{s^{2} t^{2}} \int e^{-i \frac{|w - y|^{2}}{4t}} \psi(\frac{w}{R}) e^{-i \frac{|w - z|^{2}}{4s}} dw.
\end{equation}

\noindent Now let $q(s,t,y,z) = \frac{sy + zt}{s + t} \cdot \frac{(s + t)^{1/2}}{(ts)^{1/2}}$. After making a change of variables in $w$,

\begin{equation}\label{4.12}
|K(s,t; y, z)| = \frac{C}{(s + t)^{2}} \int e^{-i |w - q(s,t,y,z)|^{2}} \psi(\frac{w}{R} \cdot \frac{(st)^{1/2}}{(s + t)^{1/2}}) dy.
\end{equation}

\noindent When $R \cdot \frac{(s + t)^{1/2}}{(st)^{1/2}} \leq 1$, H{\"o}lder's inequality implies that $|K(s,t;y,z)| \lesssim \frac{1}{(t + s)^{2}}$. For $R_{0} = R \cdot \frac{(s + t)^{1/2}}{(st)^{1/2}} > 1$, stationary phase calculations imply that for $\chi \in C_{0}^{\infty}$, $\chi = 1$ on $|x| \leq 1$, for any $N$,

\begin{equation}\label{4.13}
\int e^{-i|w - q|^{2}} (1 - \chi)(w - q) \psi(\frac{w}{R_{0}}) dw = \int ((\frac{i(w - q) \cdot \nabla}{|w - q|^{2}})^{N} e^{-i|w - q|^{2}}) (1 - \chi)(w - q) \psi(\frac{w}{R_{0}}) dw.
\end{equation}

\noindent Integrating by parts, for $N = 5$,

\begin{equation}\label{4.14}
|K(s,t; x, z)| \lesssim \frac{1}{(t + s)^{2}} \int_{|y| > 1} \frac{1}{|y|^{5}} dy \lesssim \frac{1}{(t + s)^{2}}.
\end{equation}

\noindent Therefore,







\begin{equation}\label{4.15}
\aligned
 \int_{1 < t - s_{-} < R^{2}} \int_{1 < s_{+} - t < R^{2}} \langle e^{i(t - s_{-}) \Delta} F(u_{h})(s_{-}), e^{i(t - s_{+}) \Delta} F(u_{h})(s_{+}) \rangle_{x} ds_{-} ds_{+} \\ \lesssim \int_{1 < t - s_{-} < R^{2}} \int_{1 < s_{+} - t < R^{2}} \frac{1}{(s_{+} - s_{-})^{2}} \| F(s_{-}) \|_{L_{x}^{1}} \| F(s_{+}) \|_{L_{x}^{1}} ds_{-} ds_{+} \\ \lesssim \sum_{0 \leq j \leq k \leq \ln(R)} 2^{-2k} (\int_{t - s_{+} \sim 2^{k}} \| F(u_{h})(s_{+}) \|_{L_{x}^{1}(\mathbf{R}^{4})} ds_{+}) \\ \times (\int_{t - s_{-} \sim 2^{j}} \| F(u_{h})(s_{-}) \|_{L_{x}^{1}(\mathbf{R}^{4})} ds_{-}) \lesssim \ln(R) \mathcal M(\| F(u_{h}) \|_{L_{x}^{1}(\mathbf{R}^{4})})(t)^{2}.
\endaligned
\end{equation}

\noindent Therefore, 

\begin{equation}
\langle u_{h}(t), u_{h}(t) \rangle_{x} \lesssim a(t)^{2} R^{2} + (\ln(R) + 1) \mathcal M(\| F(u_{h}) \|_{L_{x}^{1}(\mathbf{R}^{4})})(t)^{2} + 1,
\end{equation}

\noindent and since $\int [|\nabla u(t,x)|^{2} + |u(t,x)|^{4} dx] < \| W \|_{\dot{H}^{1}}^{2} + \| W \|_{L_{x}^{4}}^{4}$,

\begin{equation}
\int_{I} \int \int \psi(\frac{x - y}{R}) |u_{h}(t,y)|^{2} [|\nabla u(t,x)|^{2} + |u(t,x)|^{4}] dx dy dt \lesssim (\ln(R) + 1) K.
\end{equation}

\noindent Now notice that by the Sobolev embedding and theorem $\ref{t3.2}$,

\begin{equation}\label{lowfrequency}
\| \nabla u_{l}^{2} \|_{L_{t}^{3/2} L_{x}^{12/5}(I \times \mathbf{R}^{4})} \lesssim \| \nabla u_{l} \|_{L_{t}^{2} L_{x}^{4}(I \times \mathbf{R}^{4})}^{4/3} \| u_{l} \|_{L_{t}^{\infty} L_{x}^{4}(I \times \mathbf{R}^{4})}^{2/3} \lesssim 1,
\end{equation}

\noindent so

\begin{equation}\label{4.16}
\int_{I} \int \int |u_{l}(t,y)|^{2} \psi(\frac{(x - y)}{R}) [|\nabla u(t,x)|^{2} + |u(t,x)|^{4}] dx dy \lesssim K^{1/3} R^{10/3}.
\end{equation}

\noindent Therefore the proof of lemma $\ref{l4.1}$ is complete. $\Box$\vspace{5mm}

\noindent Now we are ready to exclude the soliton scenario. 

\begin{theorem}\label{t4.2}
Suppose $u$ is an almost periodic solution to $(\ref{equation})$ with $N(t) \equiv 1$ on $\mathbf{R}$ and with $\dot{H}^{1}$ norm below the threshold. Then $u \equiv 0$.
\end{theorem}

\noindent We prove this by constructing an interaction Morawetz estimate suited to the focusing problem. This Morawetz estimate is in the same vein as \cite{D5} and \cite{D9}.\vspace{5mm}

\noindent \emph{Proof:} Define a function $\psi \in C_{0}^{\infty}(\mathbf{R})$, $\psi$ even, $\psi = 1$ for $|x| \leq 1$ and $\psi = 0$ for $|x| > 2$. Then let

\begin{equation}\label{4.17}
\phi(x - y) = \int \psi^{2}(x - s) \psi^{2}(y - s) ds.
\end{equation}



\noindent Notice that $\phi$ is supported on $|x| \leq 4$. Then define the interaction Morawetz potential

\begin{equation}\label{4.18}
M_{R}(t) = \int |u(t,y)|^{2} \phi(\frac{x - y}{R}) (x - y)_{j} Im[\bar{u} \partial_{j} u](t,x) dx dy.
\end{equation}

\noindent By H\"older's inequality and Sobolev embedding,

\begin{equation}\label{4.19}
\sup_{t \in I} |M_{R}(t)| \lesssim R^{4}.
\end{equation}

\noindent By direct calculation,

\begin{equation}\label{4.20}
\frac{d}{dt} M_{R}(t) = 2 \int |u(t,y)|^{2} \phi(\frac{x - y}{R}) [|\nabla u(t,x)|^{2} - |u(t,x)|^{4}] dx dy
\end{equation}

\begin{equation}\label{4.21}
-2 \int Im[\bar{u} \partial_{j} u](t,y) \phi(\frac{x - y}{R}) Im[\bar{u} \partial_{j} u](t,x) dx dy
\end{equation}

\begin{equation}\label{momentum1}
+ 2 \int |u(t,y)|^{2} \phi'(\frac{x - y}{R}) \frac{(x - y)_{j} (x - y)_{k}}{|x - y| R} [Re(\partial_{j} \bar{u} \partial_{k} u)(t,x) - \frac{1}{4} \delta_{jk} |u(t,x)|^{4}] dx dy
\end{equation}

\begin{equation}\label{momentum2}
- 2 \int Im[\bar{u} \partial_{k} u](t,y) \phi'(\frac{x - y}{R}) \frac{(x - y)_{j} (x - y)_{k}}{|x - y| R} Im[\bar{u} \partial_{j} u](t,x) dx dy
\end{equation}

\begin{equation}\label{mass}
- \frac{1}{2} \int |u(t,y)|^{2} \Delta[4 \phi(\frac{x - y}{R}) + \phi'(\frac{x - y}{R}) \frac{|x - y|}{R}] |u(t,x)|^{2} dx dy.
\end{equation}

\noindent First consider $(\ref{momentum1}) + (\ref{momentum2})$. Take $R_{0} = K^{1/5}$. By the support of $\phi(x)$,

\begin{equation}\label{4.24}
\int_{1 \leq R \leq R_{0}} \frac{1}{R} |\phi'(\frac{x - y}{R}) \frac{(x - y)_{j} (x - y)_{k}}{|x - y| R}| dR \lesssim 1,
\end{equation}

\noindent and is supported on $|x - y| \lesssim R_{0}$. Therefore, by lemma $\ref{l4.1}$ and the Cauchy - Schwartz inequality,

\begin{equation}\label{4.25}
\int_{I} (\ref{momentum1}) + (\ref{momentum2}) dt \lesssim K \ln(R_{0}).
\end{equation}

\noindent Next take $(\ref{mass})$. Because $\phi(x - y)$ is supported on $|x - y| \leq 4$,

\begin{equation}\label{4.26}
\int_{1 \leq R \leq R_{0}} \frac{1}{R} |\Delta [4 \phi(\frac{x - y}{R}) + \phi'(\frac{x - y}{R}) \frac{|x - y|}{R}]| dR \lesssim \int_{1 \leq R \leq R_{0}} \frac{1}{R^{3}} \phi(\frac{x - y}{2R}) dR \lesssim \frac{1}{1 + |x - y|^{2}},
\end{equation}

\noindent and is also supported on $|x - y| \lesssim R_{0}$. Take the mollifier $\chi \in C_{0}^{\infty}(\mathbf{R}^{4})$, $\int \chi(x) = 1$, $\chi \geq 0$, and $\chi$ supported on $|x| \leq \frac{1}{4}$. Then

\begin{equation}\label{4.27}
u_{h}(t,x) = u_{h}(t,x) - \frac{1}{R^{4}} \int \chi(\frac{x - y}{R}) u_{h}(t,y) dy + \frac{1}{R^{4}} \int \chi(\frac{x - y}{R}) u_{h}(t,y) dy.
\end{equation}

\noindent By the fundamental theorem of calculus,

\begin{equation}\label{4.28}
\aligned
u_{h}(t,x) - \frac{1}{R^{4}} \int \chi(\frac{x - y}{R}) u_{h}(t,y) dy = \frac{1}{R^{4}} \int \chi(\frac{x - y}{R}) [u_{h}(t,x) - u_{h}(t,y)] dy \\ = \frac{1}{R^{4}} \int \chi(\frac{z}{R}) \int_{0}^{1} [\nabla u_{h}(t,x + sz)] \cdot z ds dz.
\endaligned
\end{equation}

\noindent Then by $(\ref{4.28})$ and lemma $\ref{l4.1}$,

\begin{equation}\label{4.29}
\aligned
\int_{1 \leq R \leq R_{0}} \frac{1}{R^{3}} \int \int_{|x - y| \sim R} |u_{h}(t,y)|^{2} |u_{h}(t,x) - \frac{1}{R^{4}} \int \chi(\frac{x - w}{R}) u_{h}(t,w) dw|^{2} dx dy dt dR \\ \lesssim \int_{I} \int_{|x - y| \lesssim R_{0}} |u_{h}(t,y)|^{2} |\nabla u_{h}(t,x)|^{2} dx dy dt dR \lesssim K \ln(R_{0}).
\endaligned
\end{equation}

\noindent By theorem $\ref{t3.2}$, since $\hat{\chi}(\xi)$ is rapidly decreasing for $|\xi| \geq 1$ and H{\"o}lder's inequality,

\begin{equation}\label{4.30}
\aligned
\int_{1 \leq R \leq R_{0}} \frac{1}{R^{3}} \int \int_{|x - y| \sim R} |\frac{1}{R^{4}} \int \chi(\frac{y - w}{R}) u(t,w) dw|^{2} |u_{h}(t,x)|^{2} dx dy dt dR \\ \lesssim \int_{1 \leq R \leq R_{0}} ( \sum_{K^{-1/4} \leq M_{1} \leq M_{2}} \frac{1}{1 + M_{1}^{4} R^{4}} \| P_{M_{1}} u \|_{L_{t}^{2} L_{x}^{4}(I \times \mathbf{R}^{4})} \frac{1}{1 + M_{2}^{4} R^{4}} \\ \times \| P_{M_{2}} u \|_{L_{t}^{2} L_{x}^{4}(I \times \mathbf{R}^{4})} \| P_{> M_{2}} u \|_{L_{t}^{\infty} L_{x}^{8/3}(I \times \mathbf{R}^{4})}^{2}) dR \lesssim K \ln(R_{0}).
\endaligned
\end{equation}

\noindent \textbf{Remark:} In fact by Bernstein's inequality, $(\ref{6.8})$, and $N(t) \geq 1$, 

\begin{equation}
\| u_{> M_{2}} \|_{L_{t}^{\infty} L_{x}^{8/3}(I \times \mathbf{R}^{4})} \lesssim o(\frac{1}{M_{2}}),
\end{equation}

\noindent where $M_{2} o(\frac{1}{M_{2}}) \rightarrow 0$ as $M_{2} \searrow 0$. Therefore, 

\begin{equation}\label{improvement}
(\ref{4.30}) \lesssim K o(\ln(R_{0})).
\end{equation}

\noindent $(\ref{improvement})$ will be crucial later. Finally, by theorem $\ref{t6.1}$, theorem $\ref{t3.2}$, and Sobolev embedding,

\begin{equation}\label{embedding}
\| u_{l}^{2} \|_{L_{t}^{2} L_{x}^{3}(I \times \mathbf{R}^{4})} \lesssim \| \nabla u_{l} \|_{L_{t}^{2} L_{x}^{4}(I \times \mathbf{R}^{4})} \| u \|_{L_{t}^{\infty} L_{x}^{3}(I \times \mathbf{R}^{4})} \lesssim 1,
\end{equation}

\noindent so by H{\"o}lder's inequality

\begin{equation}\label{4.31}
\int_{1 \leq R \leq R_{0}} \frac{1}{R^{3}} \int_{I} \int \int_{|x - y| \sim R} |u_{l}(t,x)|^{2} |u_{l}(t,y)|^{2} dx dy dt dR \lesssim K^{1/2} \int_{1 \leq R \leq R_{0}} R dR \lesssim K^{1/2} R_{0}^{2}.
\end{equation}

\noindent Now consider $(\ref{4.20})$ and $(\ref{4.21})$. Recall that

\begin{equation}\label{4.44}
\phi(\frac{x - y}{R}) = \int \psi^{2}(\frac{x}{R} - s) \psi^{2}(\frac{y}{R} - s) ds.
\end{equation}

\noindent For each $s$, $t$ there exists a $\xi(s,t)$ such that

\begin{equation}\label{4.45}
\int \psi^{2}(\frac{x}{R} - s) Im[\bar{u} \nabla e^{ix \cdot \xi(s,t)} u](t,x) dx = 0.
\end{equation}

\noindent Moreover, the quantity

\begin{equation}\label{4.46}
\int \psi^{2}(\frac{x}{R} - s) \psi^{2}(\frac{y}{R} - s) [|\nabla u(t,x)|^{2} |u(t,y)|^{2} - Im[\bar{u} \nabla u](t,x) Im[\bar{u} \nabla u](t,y) dx dy
\end{equation}

\noindent is invariant under the Galilean transformation $u \mapsto e^{-ix \cdot \xi(s,t)} u$. Therefore, for each $t$ we can take a Galilean transform on each square $\psi^{2}(x - s) \psi^{2}(y - s)$ separately to rid ourselves of the momentum squared term. Now for a fixed $t$, $s$,

\begin{equation}\label{4.47}
\aligned
\int \psi^{2}(\frac{x}{R} - s) [|\nabla e^{-ix \cdot \xi(s,t)} u(t,x)|^{2} - |u(t,x)|^{4}] dx =  \int |u(t,x)|^{2} (\psi(\frac{x}{R} - s) \Delta \psi(\frac{x}{R} - s)) dx \\ + \int |\nabla (\psi(\frac{x}{R} - s) e^{-ix \cdot \xi(s,t)} u(t,x))|^{2} dx - |\psi(\frac{x}{R} - s)u(t,x)|^{2} |u(t,x)|^{2} dx.
\endaligned
\end{equation}

\noindent By lemma $\ref{l1.12}$, $\| u \|_{\dot{H}^{1}} < (1 - \bar{\delta}) \| W \|_{\dot{H}^{1}}$, so by $(\ref{1.22})$,

\begin{equation}\label{4.48}
\| u \|_{L^{4}(\mathbf{R}^{4})} \leq (1 - \bar{\delta}) \| W \|_{L^{4}(\mathbf{R}^{4})},
\end{equation}

\noindent and therefore since $C_{4} = \frac{1}{\| W \|_{L^{4}(\mathbf{R}^{4})}}$,

\begin{equation}\label{4.49}
\aligned
\int |\nabla (\psi(\frac{x}{R} - s) e^{-ix \cdot \xi(s,t)} u(t,x))|^{2} dx - |\psi(\frac{x}{R} - s) e^{-ix \cdot \xi(s,t)} u(t,x)|^{2} |u(t,x)|^{2} dx \\ \gtrsim \bar{\delta} \int \psi(\frac{x}{R} - s)^{2} |u(t,x)|^{4} dx.
\endaligned
\end{equation}

\noindent It remains to calculate

\begin{equation}\label{4.49.1}
\int_{1 \leq R \leq R_{0}} \frac{1}{R} \int_{I} \int \int |\psi(\frac{x}{R} - s) \Delta \psi(\frac{x}{R} - s)| |u(t,x)|^{2} |u(t,y)|^{2} \psi^{2}(\frac{y}{R} - s) dx dy ds dt dR.
\end{equation}

\noindent Now if $|\frac{x}{R} - s| \sim 1$ and $|\frac{y}{R} - s| \lesssim 1$, $|x - y| \lesssim R$. Moreover, $|\Delta \psi(\frac{x}{R} - s)| \lesssim \frac{1}{R^{2}}$. Now,

\begin{equation}\label{4.50}
\aligned
\int_{1 \leq R \leq R_{0}} \frac{1}{R^{3}} \int \int_{|x - s| \sim R} \psi^{2}(y - s) |u(t,x)|^{2} |u(t,y)|^{2} dx dy ds dt dR \\ \lesssim \int_{I} \int_{|x - y| \lesssim R_{0}} \frac{1}{1 + |x - y|^{2}} |u(t,y)|^{2} |u(t,x)|^{2} dx dy dt.
\endaligned
\end{equation}

\noindent Then by $(\ref{4.29})$, $(\ref{4.30})$, and $(\ref{4.31})$,

\begin{equation}\label{4.51}
(\ref{4.50}) \lesssim K \ln(R_{0}) + K^{1/2} R_{0}^{2}.
\end{equation}



\noindent Therefore, by $(\ref{4.49})$, the fundamental theorem of calculus, $(\ref{4.19})$, $(\ref{4.25})$, $(\ref{4.29})$, $(\ref{4.30})$, $(\ref{4.31})$, and $(\ref{4.51})$,

\begin{equation}\label{4.53}
\aligned
\bar{\delta} \ln(R_{0}) \int_{I} \int_{|x - y| \leq R_{0}^{1/2}} |u(t,x)|^{4} |u(t,y)|^{2} dx dy dt - O(K \ln(R_{0})) - O(K^{1/2} R_{0}^{2}) \\  \lesssim \int_{1 \leq R \leq R_{0}} \int_{I} \frac{1}{R} \frac{d}{dt} M_{R}(t) dt dR \lesssim R_{0}^{4}.
\endaligned
\end{equation}




\noindent Therefore, if $R_{0} \leq K^{1/5}$,

\begin{equation}\label{4.54}
\bar{\delta} \int_{I} \int_{|x - y| \leq R_{0}^{1/2}} |u(t,x)|^{4} |u(t,y)|^{2} dx dt \lesssim K.
\end{equation}

\noindent Now notice that $(\ref{4.54})$ represents a logarithmic improvement over the result of lemma $\ref{l4.1}$ for the term involving $|u(t,x)|^{4}$. Moreover, because $u(t)$ lies in a compact subset of $\dot{H}^{1}$ modulo scaling symmetries, we can make a point set topology argument to prove that $\| u(t) \|_{L_{x}^{4}}$ is uniformly bounded below. Then by $(\ref{6.8})$, $(\ref{1.20})$,

\begin{equation}\label{inverse}
\int_{|x - x(t)| \leq \frac{C(\eta)}{N(t)}} |\nabla u(t,x)|^{2} dx \sim \int_{|x - x(t)| \leq \frac{C(\eta)}{N(t)}} |u(t,x)|^{4} dx.
\end{equation}

\noindent This gives a type of inverse Sobolev embedding (see \cite{CKSTT4}) that is useful to control the kinetic energy term. By $(\ref{inverse})$,

\begin{equation}\label{4.59}
\aligned
\int_{1 \leq R \leq R_{0}^{1/2}} \frac{1}{R} \int_{I} \int \int_{|y - x(t)| \leq R - C(\eta)} |u(t,y)|^{2} \psi(\frac{x - y}{R}) \frac{|x - y|}{R} |\nabla u(t,x)|^{2} dx dy dt dR \\ \lesssim \int_{I} \int \int \psi(\frac{x - y}{R_{0}}) |u(t,y)|^{2} |u(t,x)|^{4} dx dy dt.
\endaligned
\end{equation}

\noindent If $|x(t) - y| > R + C(\eta)$, then by lemma $\ref{l4.1}$ and $(\ref{6.8})$,

\begin{equation}\label{4.60}
\int_{1 \leq R \leq R_{0}} \frac{1}{R} \int \int_{|x(t) - y| > R + C(\eta)} |u(t,y)|^{2} \psi(\frac{x - y}{R})\frac{|x - y|}{R} |\nabla u(t,x)|^{2} dx dy dt dR \lesssim \eta K \ln(R_{0}).
\end{equation}

\noindent Finally, by lemma $\ref{l4.1}$,

\begin{equation}\label{4.61}
\aligned
\int_{1 \leq R \leq R_{0}^{1/2}} \frac{1}{R} \int_{I} \int_{R - C(\eta) < |x(t) - y| < R + C(\eta)} \int |\nabla u(t,x)|^{2} \psi(\frac{x - y}{R}) \frac{|x - y|}{R} |u(t,y)|^{2} dx dy dt dR \\
\lesssim \int_{I} \int \int_{|x - y| \leq 2R_{0}} \frac{C(\eta)}{C(\eta) + |x - y|} |\nabla u(t,x)|^{2} |u(t,y)|^{2} dx dy dt \lesssim K \ln(C(\eta)).
\endaligned
\end{equation}

\noindent Therefore,

\begin{equation}\label{4.62}
\aligned
\int_{I} \int_{|x - y| \leq \frac{R_{0}^{1/2}}{2}} |\nabla u(t,x)|^{2} |u(t,y)|^{2} dx dy dt \\ \lesssim \int_{1 \leq R \leq R_{0}^{1/2}} \frac{1}{R} \int \int |u(t,y)|^{2} |\nabla u(t,x)|^{2} \psi(\frac{x - y}{R}) \frac{|x - y|}{R} dx dy dt dR \\ \lesssim \eta K \ln(R_{0}) + K \ln(C(\eta)).
\endaligned
\end{equation}

\noindent Therefore,

\begin{equation}
\int_{I} (\ref{momentum1}) + (\ref{momentum2}) dt \lesssim \eta K \ln(R_{0}) + K \ln(C(\eta)) + K^{1/2} R_{0}^{2}.
\end{equation}

\noindent Next, by $(\ref{4.29})$, $(\ref{4.30})$, and $(\ref{4.31})$,

\begin{equation}\label{4.63}
\int_{1 \leq R \leq R_{0}^{1/2}} \frac{1}{R} \int_{I} (\ref{mass}) dt dR \lesssim o(\ln(R_{0})) K + \eta K \ln(R_{0}) + K \ln(C(\eta)).
\end{equation}

\noindent Then by $(\ref{4.51})$ and the fundamental theorem of calculus

\begin{equation}\label{4.64}
\aligned
\int_{I} \int_{|x - y| \leq R_{0}^{1/4}} \ln(R_{0}) |u(t,x)|^{4} |u(t,y)|^{2} dx dy dt \\ \lesssim o(\ln(R_{0})) K + \eta K \ln(R_{0}) + K \ln(C(\eta)) + R_{0}^{2} K^{1/2} + R_{0}^{4}.
\endaligned
\end{equation}

\noindent However, by $(\ref{3.2})$ this implies that either there exists a sequence $t_{n} \in \mathbf{R}$ such that $R_{0,n} \nearrow \infty$ and either

\begin{equation}\label{4.65}
\int_{|x - x(t_{n})| \leq R_{0,n}^{1/4}} |u(t_{n},x)|^{2} dx \rightarrow 0,
\end{equation}

\noindent or

\begin{equation}\label{4.66}
\int_{|x - x(t_{n})| \leq R_{0,n}^{1/4}} |u(t_{n},x)|^{4} dx \rightarrow 0.
\end{equation}

\noindent In either case, this implies that $u \equiv 0$. $\Box$



\section{Variable $N(t)$}
Now we turn to the case when $N(t)$ is free to vary.  In this case we may wish to try

\begin{equation}\label{5.1}
M(t) = \int |u(t,y)|^{2} \phi(\frac{(x - y)N(t)}{R})(x - y)_{j} Im[\bar{u} \partial_{j} u](t,x) dx dy.
\end{equation}

\noindent Everything would then proceed exactly as in theorem $\ref{t4.2}$, except that we have one additional term,

\begin{equation}\label{5.2}
\int |u(t,y)|^{2} \phi(\frac{(x - y) N(t)}{R}) \frac{|x - y| (x - y)_{j}}{R} N'(t) Im[\bar{u} \partial_{j} u](t,x) dx dy.
\end{equation}

\noindent Notice that by H\"older's inequality and $\| u \|_{L_{t}^{\infty} L_{x}^{4}(I \times \mathbf{R}^{4})} \lesssim 1$, $(\ref{5.2}) \lesssim R^{4} \frac{N'(t)}{N(t)^{5}}$. In the case that $\int \frac{N'(t)}{N(t)^{5}} dt << K$, we would be done. So this would rule out not only the case when $N(t) \equiv 1$, but also the case when $N(t)$ is a monotone function. However, $N(t)$ may be highly oscillatory. In that case, it is useful to replace $N(t)$ with $\tilde{N}(t)$, that satisfies the following conditions:

\begin{enumerate}
\item $\tilde{N}(t) \gtrsim 1$.

\item $|\tilde{N}'(t)| \lesssim \tilde{N}(t)^{3}$.

\item

\begin{equation}\label{5.3}
\int_{I} \frac{1}{\tilde{N}(t)^{2}} dt \lesssim K,
\end{equation}

\noindent and

\item 

\begin{equation}\label{5.4}
\int_{I} \frac{|\tilde{N}'(t)|}{\tilde{N}(t)^{5}} dt << K.
\end{equation}
\end{enumerate}

\noindent $\tilde{N}(t)$ will be inductively defined, using a procedure very similar to the construction in \cite{D5}. We begin with $\tilde{N}_{0}(t)$, although to simplify notation we will simply write $N_{0}(t)$.

\begin{definition}\label{d5.1}
Let

\begin{equation}\label{5.5}
\frac{1}{N_{0}(t)} = \| u_{h}(t) \|_{L_{x}^{3}(\mathbf{R}^{4})}^{3}.
\end{equation}
\end{definition}

\begin{lemma}\label{l5.2}
Possibly after modifying $N_{0}(t)$ by some function $\alpha(t)$, $N_{0}(t) \mapsto \alpha(t) N_{0}(t)$,

\begin{equation}\label{5.6}
\epsilon < \alpha(t) < \frac{1}{\epsilon},
\end{equation}

\begin{enumerate}
\item $N_{0}(t) \gtrsim 1$.

\item $|N_{0}'(t)| \lesssim N_{0}(t)^{3}$,

\noindent and

\item

\begin{equation}\label{5.7}
\int_{I} \frac{1}{N_{0}(t)^{2}} dt \lesssim K.
\end{equation}
\end{enumerate}
\end{lemma}

\noindent \emph{Proof:} $N_{0}(t) \gtrsim 1$ follows directly from theorem $\ref{t6.1}$. Notice also that by $(\ref{6.8})$ and interpolation $N_{0}(t) \lesssim N(t)$. Next, take $t_{0} \in \mathbf{R}$ and choose $N_{0} \sim N(t_{0})$ such that by Bernstein's inequality,

\begin{equation}
\| P_{> N_{0}} u(t) \|_{L_{t}^{\infty} L_{x}^{3}(\mathbf{R}^{4})} \leq 10^{-6} N(t_{0})^{-1}.
\end{equation}

\noindent By Sobolev embedding,

\begin{equation}\label{5.8}
\aligned
\frac{d}{dt} (\int |P_{\leq N_{0}} u_{h}(t,x)|^{3} dx) = (\int |P_{\leq N_{0}} u_{h}(t,x)| Re((i \Delta P_{\leq N_{0}} u_{h} + i P_{\leq N_{0}} P_{h} F(u)) P_{\leq N_{0}} \bar{u}_{h} dx) \\ \lesssim (\int |\nabla P_{\leq N_{0}} u_{h}(t,x)|^{2} |P_{\leq N_{0}} u_{h}(t,x)| dx + \int |P_{\leq N_{0}} u_{h}(t,x)|^{2} |P_{\leq N_{0}} P_{h} F(u)(t,x)| dx) \lesssim N_{0}.
\endaligned
\end{equation}

\noindent Then for $c > 0$ sufficiently small, for any $|t - t_{0}| \leq c N(t_{0})^{-2}$,

\begin{equation}\label{5.9}
\| P_{\leq N_{0}} u_{h}(t) \|_{L_{x}^{3}(\mathbf{R}^{4})} \gtrsim N(t_{0})^{-1}.
\end{equation}

\noindent Therefore, by Bernstein's inequality, for $|t_{0} - t_{1}| \leq c N(t_{0})^{-2}$,

\begin{equation}\label{5.10}
N_{0}(t_{0}) \sim N_{0}(t_{1}),
\end{equation}

\noindent and thus $|N_{0}'(t)| \lesssim N_{0}(t)^{3}$, possibly after modifying $N_{0}(t)$ by a constant.\vspace{5mm}

\noindent Finally, by theorem $\ref{t3.2}$, $(\ref{5.7})$ holds. $\Box$

\begin{theorem}\label{t5.3}
If $u$ is an almost periodic solution to $(\ref{equation})$ with $\int_{\mathbf{R}} N(t)^{-2} dt = \infty$, then $u \equiv 0$.
\end{theorem}

\noindent \emph{Proof:} Analogously to lemma $\ref{l4.1}$ define the inner product

\begin{equation}\label{5.14}
\langle f, g \rangle_{x} = \int \psi(\frac{(x - y) N_{0}(t)}{R}) f(y) \overline{g(y)} dy.
\end{equation}

\begin{equation}\label{5.15}
\int_{|t - s| > \frac{R^{4}}{N_{0}(t)^{2}}} \| e^{i(t - s) \Delta} F(u_{h})(s) ds \|_{L_{x}^{\infty}} ds \lesssim \frac{N_{0}(t)^{2}}{R^{2}} \mathcal M(\| u_{h} \|_{L^{3}}^{3})(t).
\end{equation}

\begin{equation}\label{5.16}
\| \int_{|t - s| < \frac{1}{N_{0}(t)^{2}}} e^{i(t - s) \Delta} F(u_{h})(s) ds \|_{L^{2}(\mathbf{R}^{4})} \lesssim \frac{1}{N_{0}(t)}.
\end{equation}

\noindent Next, by $(\ref{4.14})$,

\begin{equation}\label{5.17}
\int_{\frac{-R^{2}}{N_{0}(t)^{2}} + t}^{\frac{-1}{N_{0}(t)^{2}} + t} \int_{\frac{1}{N_{0}(t)^{2}} + t}^{\frac{R^{2}}{N_{0}(t)^{2}} + t} \langle e^{i(t - s_{+}) \Delta} F(u_{h})(s_{+}), e^{i(t - s_{-}) \Delta} F(u_{h})(s_{-}) \rangle_{x} ds_{+} ds_{-} \lesssim \mathcal M(\| u_{h} \|_{L_{x}^{3}}^{3})(t)^{2}.
\end{equation}

\noindent Therefore by $(\ref{4.6})$, $N_{0}(t) \gtrsim 1$, $(\ref{5.7})$, $(\ref{5.15})$, $(\ref{5.16})$, and $(\ref{5.17})$, and $(\ref{lowfrequency})$,

\begin{equation}\label{5.18}
\aligned
 \int_{I} \int |u(t,y)|^{2} \psi(\frac{(x - y) N_{0}(t)}{R}) [|\nabla u(t,x)|^{2} + |u(t,x)|^{4}] dx dy dt \\ \lesssim K \ln(R) + K^{1/2} R^{2} + (\int_{I} \frac{R^{10}}{N_{0}(t)^{10}} dt)^{1/3} \lesssim K \ln(R) + K^{1/2} R^{2} + K^{1/3} R^{10/3}.
\endaligned
\end{equation}

\noindent Next, by a calculation similar to $(\ref{4.29})$,

\begin{equation}\label{5.19}
\aligned
\int_{1 \leq R \leq R_{0}} \frac{N_{0}(t)^{2}}{R^{3}} \int_{I} \int_{|x - y| \leq 2\frac{R}{N_{0}(t)}} |u_{h}(t,y)|^{2} \\ \times |u_{h}(t,x) - \frac{N_{0}(t)^{4}}{R^{4}} \int \chi(\frac{(x - w) N_{0}(t)}{R}) u_{h}(t,w) dw|^{2} dx dy dt dR \\ \lesssim \int_{I} \int_{|x - y| \leq 2 \frac{R_{0}}{N_{0}(t)}} |\nabla u_{h}(t,x)|^{2} |u_{h}(t,y)|^{2} dx dy dt \lesssim \ln(R_{0}) K + R_{0}^{2} K^{1/2}.
\endaligned
\end{equation}

\noindent Also,

\begin{equation}\label{5.20}
\aligned
\int_{1 \leq R \leq R_{0}} \frac{N_{0}(t)^{2}}{R^{3}} \int_{I} \int_{|x - y| \leq 2\frac{R}{N_{0}(t)}} |\frac{N_{0}(t)^{4}}{R^{4}} \int \chi(\frac{(y - w) N_{0}(t)}{R}) u_{h}(t,w) dw|^{2} \\ \times |u_{h}(t,x) - \frac{N_{0}(t)^{4}}{R^{4}} \int \chi(\frac{(x - w) N_{0}(t)}{R}) u_{h}(t,w) dw|^{2} dx dy dt dR \\ \lesssim \int_{I} \int_{|x - y| \leq 2 \frac{R_{0}}{N_{0}(t)}} |\nabla u_{h}(t,x)|^{2} |u_{h}(t,y)|^{2} dx dy dt \lesssim \ln(R_{0}) K + R_{0}^{2} K^{1/2}.
\endaligned
\end{equation}

\noindent Finally, by H{\"o}lder's inequality, the fact that $\hat{\chi}$ is rapidly decreasing for $|\xi| \geq 1$, $(\ref{6.8})$, and theorem $\ref{t3.2}$,

\begin{equation}\label{5.21}
\aligned
\int_{I} \frac{N_{0}(t)^{2}}{R^{3}} \int_{|x - y| \leq \frac{2R}{N_{0}(t)}} |\frac{N_{0}(t)^{4}}{R^{4}} \int \chi(\frac{N_{0}(t)(y - w)}{R}) u_{h}(t,w) dw|^{2} \\ \times |\frac{N_{0}(t)^{4}}{R^{4}} \int \chi(\frac{(x - w) N_{0}(t)}{R}) u_{h}(t,w) dw|^{2} dx dy dt \\ \lesssim \int_{I} \sum_{K_{1}^{-1/4} \leq N_{1} \leq N_{2} \leq N_{3} \leq N_{4}} \frac{R N_{0}(t)^{2}}{N_{0}(t)^{4} + R^{4} N_{4}^{4}}  \| P_{N_{1}} u_{h}(t) \|_{L_{x}^{\infty}(\mathbf{R}^{4})} \| P_{N_{2}} u_{h}(t) \|_{L_{x}^{\infty}(\mathbf{R}^{4})} \\ \times \| P_{N_{3}} u_{h}(t) \|_{L_{x}^{2}(\mathbf{R}^{4})} \| P_{N_{4}} u_{h}(t) \|_{L_{x}^{2}(\mathbf{R}^{4})} dt \\ \lesssim o(\frac{1}{R}) \| \sup_{M \geq K_{1}^{-1/4}} M^{-2} \| P_{M} u_{h}(t) \|_{L_{x}^{\infty}(\mathbf{R}^{4})} \|_{L_{t}^{2}(I)}^{2} \| \nabla u_{h} \|_{L_{t}^{\infty} L_{x}^{2}(I \times \mathbf{R}^{4})}^{2} \lesssim K o(\frac{1}{R}).
\endaligned
\end{equation}

\begin{equation}\label{5.22}
\int_{1 \leq R \leq R_{0}} o(\frac{1}{R}) K dR \lesssim K o(\ln(R_{0})).
\end{equation}


\noindent Finally, by H\"older's inequality, Sobolev embedding, theorem $\ref{t3.2}$, theorem $\ref{t6.1}$, and $(\ref{embedding})$,

\begin{equation}\label{5.23}
\aligned
\int_{1 \leq R \leq R_{0}} \frac{N_{0}(t)^{2}}{R^{3}} \int_{I} \int_{|x - y| \lesssim \frac{R}{N_{0}(t)}} |u_{l}(t,y)|^{2} |u_{l}(t,x)|^{2} dx dy dt dR \\ \lesssim \int_{1 \leq R \leq R_{0}} (\int_{I} \frac{R^{2}}{N_{0}(t)^{4}} dt)^{1/2} dR \lesssim R_{0}^{2} K^{1/2}.
\endaligned
\end{equation}

\noindent Therefore, by $(\ref{5.19})$ - $(\ref{5.23})$,

\begin{equation}\label{5.24}
\int_{1 \leq R \leq R_{0}} \int_{I} \frac{N_{0}(t)^{2}}{R^{3}} \int_{|x - y| \leq \frac{2R}{N_{0}(t)}} |u(t,y)|^{2} |u(t,x)|^{2} dx dy dt dR \lesssim K \ln(R_{0}) + K^{1/2} R_{0}^{2}.
\end{equation}

\noindent Now we will define an interaction Morawetz estimate with $N_{m}(t) \geq N_{0}(t)$, $N_{m}(t)$ becoming progressively smoother in time after each iteration. Let

\begin{equation}\label{5.25}
M_{R}(t) = \int_{I} \int \int \psi(\frac{(x - y) N_{m}(t)}{R}) (x - y)_{j} |u(t,y)|^{2} Im[\bar{u} \partial_{j} u](t,x) dx dy.
\end{equation}

\noindent Then $|M_{R}(t)| \lesssim \frac{R^{4}}{N_{m}(t)^{4}} \lesssim R^{4}$ and

\begin{equation}\label{5.26}
\frac{d}{dt} M_{R}(t) = 2 \int \int \psi(\frac{(x - y) N_{m}(t)}{R}) |u(t,y)|^{2} [|\nabla u(t,x)|^{2} - |u(t,x)|^{4}] dx dy
\end{equation}

\begin{equation}\label{5.27}
- 2 \int \int \psi(\frac{(x - y) N_{m}(t)}{R}) Im[\bar{u} \partial_{j} u](t,x) Im[\bar{u} \partial_{j} u](t,y) dx dy
\end{equation}

\begin{equation}\label{5.28}
+ 2 \int \int \psi'(\frac{(x - y) N_{m}(t)}{R}) \frac{(x - y)_{j} (x - y)_{k}}{|x - y| R} |u(t,y)|^{2} [Re(\partial_{j} \bar{u} \partial_{k} u)(t,x) - \frac{\delta_{jk}}{4} |u(t,x)|^{4}] dx dy
\end{equation}

\begin{equation}\label{5.28.1}
- 2 \int \int Im[\bar{u} \partial_{k} u](t,y) \psi'(\frac{(x - y) N_{m}(t)}{R}) \frac{(x - y)_{j} (x - y)_{k}}{|x - y| R} Im[\bar{u} \partial_{j} u](t,x) dx dy
\end{equation}

\begin{equation}\label{5.29}
+ \frac{1}{2} \int \int |u(t,y)|^{2} \Delta [4 \psi(\frac{(x - y) N_{m}(t)}{R}) + \psi'(\frac{(x - y) N_{m}(t)}{R}) \frac{|x - y| N_{m}(t)}{R}] |u(t,x)|^{2} dx dy
\end{equation}

\begin{equation}\label{5.30}
+ \int \int \psi'(\frac{(x - y) N_{m}(t)}{R}) \frac{(x - y)_{j} |x - y| N_{m}'(t)}{R} |u(t,y)|^{2} Im[\bar{u} \partial_{j} u](t,x) dx dy.
\end{equation}

\noindent Now by $(\ref{4.47}) - (\ref{4.49})$,

\begin{equation}\label{5.31}
\aligned
(\ref{5.26}) + (\ref{5.27}) \geq \bar{\delta} \int \psi(\frac{(x - y) N_{m}(t)}{R}) |u(t,y)|^{2} |u(t,x)|^{4} dx dy \\ - \frac{C N_{m}(t)^{2}}{R^{2}} \int_{|x - y| \leq 2\frac{R}{N_{m}(t)}} |u(t,y)|^{2} |u(t,x)|^{2} dx dy.
\endaligned
\end{equation}

\noindent Therefore, by $(\ref{5.18})$ - $(\ref{5.23})$, for $R_{0} \leq K^{1/5}$,

\begin{equation}\label{5.32}
\aligned
\int_{I} \ln(R_{0}) \int \int_{|x - y| \leq \frac{R_{0}^{11/12}}{N_{m}(t)}} |u(t,y)|^{2} |u(t,x)|^{4} dx dy dt - K \ln(R_{0}) - K^{1/2} R_{0}^{2}  \\ - \int_{1 \leq R \leq R_{0}} \frac{1}{R} \int_{I} \int \psi'(\frac{(x - y) N_{m}(t)}{R}) \frac{(x - y)_{j} |x - y| N_{m}'(t)}{R} \\ \times |u(t,y)|^{2} Im[\bar{u} \partial_{j} u](t,x) dx dy dt dR \\ \gtrsim \int_{I} \int_{1 \leq R \leq R_{0}} \frac{1}{R} \frac{d}{dt} M_{R}(t) dt \lesssim R_{0}^{4}.
\endaligned
\end{equation}

\noindent Now we apply a smoothing procedure tto make $N_{m}(t)$ much smoother than $N_{0}(t)$. See \cite{D5} for a similar procedure. Partition $I$ into subintervals $J_{k}$ such that $\int_{J_{k}} N_{0}(t)^{2} dt = c$ for some $c << 1$. Then let 

\begin{equation}
N(J_{k}) = \sup \{ 2^{j} : j \in \mathbf{Z}, \hspace{5mm} 2^{j} \leq N(t) \hspace{5mm} \forall t \in J_{k} \}.
\end{equation}

\noindent For $c$ sufficiently small, if $J_{k}$ and $J_{k + 1}$ are adjacent intervals then $|N_{0}'(t)| \lesssim N_{0}^{3}(t)$ implies

\begin{equation}\label{5.33}
\frac{N(J_{k})}{N(J_{k + 1})} = 1, \hspace{5mm} \frac{1}{2}, \hspace{5mm} \text{ or } 2
\end{equation}

\noindent and $|J_{k}| \sim N(J_{k})^{-2}$. Then choose $N_{1}(t_{k}) = N(J_{k})$, where $t_{k}$ is the midpoint of $J_{k}$ and let $N_{1}(t)$ be the linear interpolation between these midpoints. Then

\begin{equation}\label{5.34}
\int_{I} \frac{|N_{1}'(t)|}{N_{1}(t)^{5}} dt \lesssim K.
\end{equation}

\noindent Now we iteratively obtain $N_{l + 1}(t)$ from $N_{l}(t)$ using the smoothing algorithm.

\begin{definition}[Smoothing algorithm]\label{d5.3}
An interval $J_{k}$ is called upward sloping if $\frac{N(J_{k})}{N(J_{k + 1})} = \frac{1}{2}$, downward sloping if $\frac{N(J_{k})}{N(J_{k + 1})} = 2$, and flat if $\frac{N(J_{k})}{N(J_{k + 1})} = 1$. We call $J$ a valley if $J = J_{l} \cup J_{l + 1} \cup ... \cup J_{l + m}$, $J_{l}$ is downward sloping, $J_{l + m}$ is upward sloping, and $J_{l + 1}$, ..., $J_{l + m - 1}$ are constant intervals. We call $J$ a peak if $J = J_{l} \cup J_{l + 1} \cup ... \cup J_{l + m}$, $J_{l}$ is upward sloping, $J_{l + m}$ is downward sloping, and $J_{l + 1}$, ..., $J_{l + m - 1}$ are constant intervals.\vspace{5mm}

\noindent \textbf{Remark:} $N_{1}(t)$ is monotone in between consecutive peaks and valleys. Moreover, we cannot have two peaks without a valley in between, or two valleys without a peak in between.\vspace{5mm}

\noindent Now if

\begin{equation}\label{5.35}
J = J_{l} \cup ... \cup J_{m + 1},
\end{equation}

\noindent is a valley let $N_{2}(t) = N_{1}(t_{l}) = N_{j}(t_{l})$ for all $t_{l} < t < t_{l + m}$. Otherwise let $N_{2}(t) = N_{1}(t)$.
\end{definition}

\noindent Likewise construct $N_{3}(t)$ using the above algorithm with $N_{1}(t)$ replaced by $N_{2}(t)$. Now, by the fundamental theorem of calculus, if $N_{j}(t)$ is monotone on an interval $J_{0}$,

\begin{equation}\label{5.36}
\int_{J_{0}} \frac{|N_{j}'(t)|}{N_{j}(t)^{5}} dt \leq (\inf_{t \in J_{0}} N_{j}(t))^{-4}.
\end{equation}

\noindent Therefore, by induction, $(\ref{5.36})$, the smoothing algorithm, and the fundamental theorem of calculus,

\begin{equation}\label{5.37}
\int_{I} \frac{|N_{m}'(t)|}{N_{m}^{5}(t)} dt \leq 2^{-4m + 4} \int_{I} \frac{|N_{1}'(t)|}{N_{1}(t)^{5}} dt + 2.
\end{equation}

\noindent Observe also that 

\begin{equation}\label{5.38}
N_{1}(t) \leq N_{m}(t) \leq 2^{m - 1} N_{1}(t).
\end{equation}

\noindent Next, observe that by the definition of $N_{0}(t)$, $N_{1}(t) \sim N_{0}(t)$, H{\"o}lder's inequality, and the fact that as in \cite{D5}, either $N_{m}(t) = N_{0}(t)$ or $N_{m}'(t) = 0$,

\begin{equation}\label{5.39}
\aligned
\int_{I} \int_{|x - y| \leq \frac{R}{N_{m}(t)}} |u(t,y)|^{2} \frac{|x - y|^{2} |N_{m}'(t)|}{R} |\nabla u(t,x)| |u(t,x)| dx dy dt \\ \lesssim R^{3} \int_{I} \frac{|N_{m}'(t)|}{N_{m}(t)^{4}} \| u_{h}(t) \|_{L_{x}^{3}(\mathbf{R}^{4})}^{3} \| \nabla u \|_{L_{x}^{2}(\mathbf{R}^{4})} dt + \int_{I} \frac{R^{5} |N_{m}'(t)|}{N_{m}(t)^{6}} \| u_{l}(t) \|_{L_{x}^{6}(\mathbf{R}^{4})}^{3} \| \nabla u(t) \|_{L_{x}^{2}(\mathbf{R}^{4})} dt \\ \lesssim \int_{I} \frac{|N_{m}'(t)|}{N_{m}(t)^{5}} R^{3} dt \lesssim 2^{-4m + 4} K R^{3} + R^{3} + K^{1/2} R^{5}.
\endaligned
\end{equation}

\noindent Choose $m$ so that $2^{-4m} = R^{-3}$. Then by $(\ref{5.32})$,

\begin{equation}\label{5.40}
\ln(R_{0}) \int_{I} \int \int_{|x - y| \leq \frac{R_{0}^{11/12}}{N_{m}(t)}} |u(t,y)|^{2} |u(t,x)|^{4} dx dy dt \lesssim K \ln(R_{0}) + K^{1/2} R_{0}^{2} + R_{0}^{4} + K + R_{0}^{3} + K^{1/2} R_{0}^{5}.
\end{equation}

\noindent Therefore, if $R_{0} \leq K^{1/10}$,

\begin{equation}\label{5.41}
\int_{I} \int \int_{|x - y| \leq \frac{R_{0}^{11/12}}{N_{m}(t)}} |u(t,y)|^{2} |u(t,x)|^{4} dx dy dt \lesssim K.
\end{equation}

\noindent Now since $N_{m}(t) \leq 2^{m - 1} N(t) \leq R_{0}^{3/4} N(t)$,

\begin{equation}\label{5.42}
\int_{I} \int \int_{|x - y| \leq \frac{R_{0}^{1/6}}{N_{1}(t)}} |u(t,y)|^{2} |u(t,x)|^{4} dx dy dt \lesssim K.
\end{equation}

\noindent Now, by the inverse Sobolev embedding, $(\ref{inverse})$, $N_{1}(t) \lesssim N(t)$, $(\ref{5.42})$ implies

\begin{equation}\label{5.43}
\int_{1 \leq R \leq R_{0}^{1/6}} \frac{1}{R} \int_{I} \int_{|x - y| \leq \frac{R}{N_{1}(t)}} \int_{|y - y(t)| \leq \frac{R - C(\eta)}{N_{1}(t)}} \frac{N_{1}(t) |x - y|}{R} |\nabla u(t,x)|^{2} |u(t,y)|^{2} dx dy dt dR \lesssim K.
\end{equation}

\noindent Also by $(\ref{5.24})$,

\begin{equation}\label{5.44}
\int_{1 \leq R \leq R_{0}^{1/6}} \frac{1}{R} \int_{I} \int_{|x - y| \leq \frac{R}{N_{1}(t)}} \int_{|y - y(t)| > \frac{R + C(\eta)}{N_{1}(t)}} \frac{N_{1}(t) |x - y|}{R} |\nabla u(t,x)|^{2} |u(t,y)|^{2} dx dy dt dR \lesssim \eta K \ln(R_{0}).
\end{equation}

\noindent Finally, by $(\ref{5.24})$,

\begin{equation}\label{5.45}
\aligned
\int_{1 \leq R \leq R_{0}^{1/6}} \frac{1}{R} \int_{I} \int_{|x - y| \leq \frac{R}{N_{1}(t)}} \int_{\frac{R - C(\eta)}{N_{1}(t)} \leq |y - y(t)| \leq \frac{R + C(\eta)}{N(t)}} \frac{N_{1}(t) |x - y|}{R} |\nabla u(t,x)|^{2} |u(t,y)|^{2} dx dy dt dR \\ \lesssim \int_{I} \int_{|x - y| \leq \frac{R_{0}^{1/6}}{N_{1}(t)}} \frac{C(\eta)}{C(\eta) + N_{1}(t) |x - y|} |\nabla u(t,x)|^{2} |u(t,y)|^{2} dx dy dt \lesssim \ln(C(\eta)) K.
\endaligned
\end{equation}

\noindent Then by taking $m$ such that $2^{4m} = R_{0}^{5/9}$, by $(\ref{5.26})$ - $(\ref{5.30})$, $(\ref{5.43})$ - $(\ref{5.45})$, since $N_{m}(t) \geq N_{1}(t)$,

\begin{equation}\label{5.46}
\aligned
R_{0}^{2/3} \gtrsim \int_{1 \leq R \leq \frac{R_{0}^{1/6}}{N_{m}(t)}} \frac{1}{R} \int_{I} \frac{d}{dt} M_{R}(t) dt dR \\ \gtrsim \ln(R_{0}) \int_{I} \int \int_{|x - y| \leq \frac{R_{0}^{11/72}}{N_{m}(t)}} |u(t,y)|^{2} |u(t,x)|^{4} dx dy dt dR \\ - \ln(C(\eta)) K - \eta K \ln(R_{0}) - R_{0}^{1/2} - R_{0}^{-1/18} K - K^{1/2} R_{0}^{5/6}.
\endaligned
\end{equation}

\noindent Also since $N_{m}(t) \leq 2^{5/36} N_{1}(t)$,

\begin{equation}\label{5.47}
\aligned
\int_{I} \int \int_{|x - y| \leq \frac{R_{0}^{1/72}}{N_{1}(t)}} |u(t,y)|^{2} |u(t,x)|^{4} dx dy dt \\ \lesssim \eta K + \frac{C(\eta)}{\ln(R_{0})} K + \frac{o(\ln(R_{0}))}{\ln(R_{0})} K + R_{0}^{-1/18} K + K^{1/2} R_{0}^{5/6} + R_{0}^{1/2} + R_{0}^{2/3}.
\endaligned
\end{equation}

\noindent Therefore, since $N_{1}(t) \leq N(t)$ there exists a sequence $t_{n} \in \mathbf{R}$, $R_{n} \nearrow \infty$ such that either

\begin{equation}\label{5.48}
N(t_{n})^{2} \int_{|x - x(t)| \leq \frac{R_{n}}{N(t)}} |u(t,x)|^{2} dx \rightarrow 0,
\end{equation}

\noindent or

\begin{equation}\label{5.49}
\int_{|x - x(t)| \leq \frac{R_{n}}{N(t)}} |u(t,x)|^{4} dx \rightarrow 0.
\end{equation}

\noindent However, if $u$ lies in a precompact set modulo scaling and translation symmetries, $\| u(t) \|_{L_{x}^{4}(\mathbf{R}^{4})} \gtrsim \| u(t) \|_{\dot{H}^{1}(\mathbf{R}^{4})}$. Therefore, $(\ref{5.48})$ or $(\ref{5.49})$ imply that $u \equiv 0$. $\Box$

\nocite*
\bibliographystyle{plain}

\begin{thebibliography}{[00]}
\bibitem{Aubin}
	\newblock T. Aubin,
	\newblock ``Equations differentielles non lineaires et probleme de Yamabe concernant la courbure scalaire",
	\newblock \textit{J. Math. Pures Appl.} (9), 55, 1976, 3, 269 -- 296.

\bibitem{BahGer}
	\newblock H. Bahouri and P. Gerard,
	\newblock ``High frequency approximations of solutions to critical nonlinear wave equations'',
	\newblock \textit{American Journal of Mathematics}, \textbf{121} no. 1 (1999) 131 -- 175.

\bibitem{BerCaz}
	\newblock H. Berestycki and T. Cazenave,
	\newblock ``Instabilit\'e des \'etats stationnaires dans les \'equations de {S}chr\"odinger et de {K}lein-{G}ordon non li\'eaires",
	\newblock \textit{Comptes Rendus des S\'eances de l'Acad\'emie des Sciences. S\'erie I. Math\'ematique}, \textbf{293} no. 9 (1981) 489 -- 492.

\bibitem{B4}
	\newblock J. Bourgain,
	\newblock ``Global well - posedness of defocusing critical nonlinear Schr{\"o}dinger equation in the radial case'',
	\newblock \textit{Journal of the American Mathematical Society}, \textbf{12} (1999) 145 -- 171.

\bibitem{B3}
	\newblock J. Bourgain,
	\newblock \textit{Global {S}olutions of {N}onlinear {S}chr{\"o}dinger {E}quations},
	\newblock American Mathematical Society, \textbf{46} American Mathematical Society Colloquium Publications, Providence, RI, 1999.

\bibitem{Bulut1} 
	\newblock A. Bulut,
	\newblock ``The defocusing energy - supercritical cubic nonlinear wave equation in dimension five",
	\newblock \textit{arxiv 1112.0629}.

\bibitem{Bulut2}
	\newblock A. Bulut,
	\newblock ``The radial defocusing energy - supercritical cubic nonlinear wave equation in $\mathbf{R}^{1 + 5}$",
	\newblock \textit{arxiv 1104.2002}.

\bibitem{Bulut3}
	\newblock A. Bulut,
	\newblock ``Global well - posedness and scattering for the defocusing energy - supercritical cubic nonlinear wave equation",
	\newblock \textit{Journal of Functional Analysis} \textbf{263} no. 6 (2012) 1609 -- 1660.

\bibitem{Caz1}
	\newblock T. Cazenave,
	\newblock \textit{Semilinear Schr{\"o}dinger equations},
	\newblock Courant Lecture Notes in Mathematics \textbf{10}, New York University, Courant Institute of Mathematical Sciences, AMS, Providence, RI, 2003.

\bibitem{CaWe}
    	\newblock T. Cazenave and F. B. Weissler,
    	\newblock ``The {C}auchy problem for the nonlinear {S}chr\"odinger equation in {$H^1$}'',
    	\newblock \textit{Manuscripta Math.}, \textbf{61} (1988) 477 -- 494.

\bibitem{CaWe1}
    	\newblock T. Cazenave and F. B. Weissler,
    	\newblock ``The {C}auchy problem for the critical nonlinear {S}chr\"odinger equation in {$H^s$}'',
    	\newblock \textit{Nonlinear Anal.}, \textbf{14} (1990) 807 -- 836.

\bibitem{CKSTT2}
	\newblock J. Colliander, M. Keel, G. Staffilani, H. Takaoka, and T. Tao,
	\newblock ``Global existence and scattering for rough solutions of a nonlinear {S}chr{\"o}dinger equation on $\mathbf{R}^{3}$'',
        \newblock \textit{Communications on Pure and Applied Mathematics}, \textbf{21} (2004) 987 - 1014.

\bibitem{CKSTT4}
	\newblock J. Colliander, M. Keel, G. Staffilani, H. Takaoka, and T. Tao,
	\newblock ``Global well - posedness and scattering for the energy - critical nonlinear {S}chr\"odinger equation on $\mathbf{R}^{3}$'',
        \newblock \textit{Annals of Mathematics. Second Series}, \textbf{167} (2008) 767 - 865.

\bibitem{D2}
    	\newblock B. Dodson,
    	\newblock ``Global well - posedness and scattering for the defocusing {$L^{2}$} - critical nonlinear {S}chr{\"o}dinger equation when $d \geq 3$'',
    	\newblock \textit{Journal of the American Mathematical Society}, \textbf{25} no. 2 (2012) 429 -- 463.

\bibitem{D3}
    	\newblock B. Dodson,
    	\newblock ``Global well-posedness and scattering for the defocusing, $L^{2}$ - critical nonlinear Schr{\"o}dinger equation when $d = 2$, preprint,
    	\newblock \textit{arXiv:1006.1375},
	
\bibitem{D4}
    	\newblock B. Dodson,
    	\newblock ``Global well-posedness and scattering for the defocusing, $L^{2}$ - critical nonlinear Schr{\"o}dinger equation when $d = 1$, preprint,
    	\newblock \textit{arXiv:1010.0040},	

\bibitem{D5}
    	\newblock B. Dodson,
    	\newblock ``Global well-posedness and scattering for the mass critical nonlinear Schr{\"o}dinger equation with mass below the mass of the ground state'', preprint,
    	\newblock \textit{arXiv:1104.1114v2},

\bibitem{D9}
	\newblock B. Dodson,
	\newblock ``Global well - posedness and scattering for the defocusing, mass - critical generalized {K}d{V} equation'', preprint,
	\newblock \textit{arxiv 1304.8025}.

\bibitem{GV}
	\newblock J. Ginibre and G. Velo,
	\newblock ``Smoothing properties and retarded estimates for some dispersive evolution equations'',
        \newblock \textit{Communications in Mathematical Physics}, \textbf{144} no. 1 (1992) 163 -- 188.

\bibitem{Glassey}
	\newblock R. T. Glassey,
	\newblock ``On the blowing up of solutions to the {C}auchy problem for nonlinear {S}ch{\"o}dinger equations",
	\newblock \textit{Journal of Mathematical Physics}, \textbf{18} no. 9 (1977) 1794 -- 1797.

\bibitem{Gril}
	\newblock M. Grillakis,
	\newblock ``On nonlinear Schr{\"o}dinger equations'',
        \newblock \textit{Communications in Partial Differential Equations}, \textbf{25} no. 9 - 10 (2000) 1827 -- 1844.

\bibitem{KT}
	\newblock M. Keel and T. Tao,
	\newblock ``Endpoint {S}trichartz Estimates''
	\newblock \textit{American Journal of Mathematics} \textbf{120} no. 4 - 6 (1998) 945 -- 957.

\bibitem{KM1}
	\newblock C. Kenig and F. Merle,
	\newblock ``Global well-posedness, scattering, and blow-up for the energy-critical, focusing nonlinear Schr{\"o}dinger equation in the radial case,''
	\newblock \textit{Inventiones Mathematicae} \textbf{166} no. 3 (2006) 645--675.

\bibitem{Keraani1}
	\newblock S. Keraani,
	\newblock ``On the defect of compactness for the {S}trichartz estimates of the {S}chr\"odinger equations'',
	\newblock \textit{Journal of Differential Equations} \textbf{175} no. 2 (2001) 353 -- 392.

\bibitem{KTV}
    	\newblock R. Killip, T. Tao, and M. Visan,
     	\newblock ``The cubic nonlinear {S}chr\"odinger equation in two dimensions with radial data",
   	\newblock \textit{Journal of the European Mathematical Society (JEMS)}, \textbf{11} no. 6 (2009) 1203 -- 1258.

\bibitem{KilVis}
	\newblock R. Killip and M. Visan,
	\newblock ``Nonlinear Schr{\"o}dinger Equations at Critical Regularity",
	\newblock \textit{Unpublished lecture notes} , Clay Lecture Notes (2009): http://www.math.ucla.edu/~visan/lecturenotes.html.

\bibitem{KV1}
	\newblock R. Killip and M. Visan,
	\newblock ``The focusing energy - critical nonlinear {S}chr{\"o}dinger equation in dimensions five and higher'',
	\newblock \textit{American Journal of Mathematics}, textbf{132}, no. 2 (2010) 361 -- 424.

\bibitem{KV4}
	\newblock R. Killip and M. Visan,
	\newblock ``The radial defocusing energy - supercritical nonlinear wave equation in all space dimensions",
	\newblock \textit{Proceedings of the American Mathematical Society} \textbf{139} no. 5 (2011) 1805 -- 1817.

\bibitem{KV3}
	\newblock R. Killip and M. Visan,
	\newblock ``The defocusing energy - supercritical nonlinear wave equation in three space dimensions",
	\newblock \textit{Transactions of the American Mathematical Society} \textbf{363} no. 7 (2011) 3893 -- 3934.

\bibitem{KV2}
	\newblock R. Killip and M. Visan,
	\newblock ``Global well - posedness and scattering for the defocusing quintic {NLS} in three dimensions",
	\newblock \textit{Analysis and PDE} \textbf{5} no. 4 (2012) 855 -- 885.

\bibitem{KVZ}
    	\newblock R. Killip, M. Visan, and X. Zhang,
     	\newblock ``The mass-critical nonlinear {S}chr\"odinger equation with radial data in dimensions three and higher",
  	\newblock \textit{Analysis \& PDE}, \textbf{1} no. 2 (2008) 229 -- 266.

\bibitem{Merle}
	\newblock F. Merle,
     	\newblock ``Determination of blow-up solutions with minimal mass for nonlinear {S}chr\"odinger equations with critical 		power",
	\newblock \textit{Duke Mathematical Journal}, \textbf{69} no. 2 (1993) 427 -- 454.

\bibitem{MMZ}
	\newblock C. Miao, J. Murphy, and J. Zheng,
	\newblock ``The defocusing energy - supercritical {NLS} in four space dimensions",
	\newblock \textit{Journal of Functional Analysis} \textbf{267} (2014) 1662 -- 1724.

\bibitem{RhV}
	\newblock E. Rhyckman and M. Visan,
	\newblock ``Global well-posedness and scattering for the defocusing energy-critical nonlinear {S}chr\"odinger equation in {$\Bbb R^{1+4}$}'',
	\newblock \textit{American Journal of Mathematics} \textbf{129}, no. 1 (2007): 1 -- 60.

\bibitem{St1} 
     	\newblock E. M. Stein,
     	\newblock \textit{Singular Integrals and Differentiability Properties of functions},
     	\newblock Princeton University Press, Princeton, NJ, 1970.

\bibitem{St} 
     	\newblock E. M. Stein,
     	\newblock \textit{Harmonic Analysis: Real-variable Methods, Orthogonality, and Oscillatory Integrals},
     	\newblock Princeton University Press, Princeton, NJ, 1993.

\bibitem{Stri}
	\newblock R. S. Strichartz,
	\newblock ``Restrictions of {F}ourier transforms to quadratic surfaces and decay of solutions of wave equations'',
	\newblock \textit{Duke Mathematical Journal} \textbf{44} no. 3 (1977) 705 - 714.

\bibitem{Talenti}
	\newblock G. Talenti,
	\newblock ``Best constant in {S}obolev inequality'',
	\newblock \textit{Annali di Matematica Pura ed Applicata. Serie Quarta} \textbf{110} (1976) 353 -- 372.

\bibitem{TaoT1} 
     	\newblock T. Tao,
     	\newblock ``On the asymptotic behavior of large radial data for a focusing nonlinear Schr{\"o}dinger equation'',
     	\newblock \textit{Dynamics of PDE}, \textbf{1} no. 1 (2004) 1 -- 48.

\bibitem{TerryTao} 
     	\newblock T. Tao,
     	\newblock ``Global well - posedness and scattering for the higher - dimensional energy - critical nonlinear {S}chr{\"o}dinger equation for radial data'',
     	\newblock \textit{New York Journal of Mathematics} \textbf{11} (2005) 57 -- 80.

\bibitem{Tao} 
     	\newblock T. Tao,
     	\newblock \textit{Nonlinear Dispersive Equations. Local and Global Analysis},
     	\newblock CBMS Regional Conference Series in Mathematics \textbf{104} Published for the Conference Board of the Mathematical Sciences, Washington, DC, 2006.

\bibitem{TVZ}
	\newblock T. Tao, M. Visan, and X. Zhang.
	\newblock ``The nonlinear {S}chr\"odinger equation with combined power-type nonlinearities'',
	\newblock \textit{Communications in Partial Differential Equations}, \textbf{32} no. 7-9 (2007) 1281--1343.

\bibitem{T3} 
     	\newblock M. E. Taylor,
     	\newblock \textit{Pseudodifferential Operators and Nonlinear PDE},
     	\newblock Birkh{\"a}user, Boston, 1991.

\bibitem{T4} 
     	\newblock M. E. Taylor,
     	\newblock \textit{Tools for PDE. Pseudodifferential operators, paradifferential Operators, and layer potentials},
     	\newblock Mathematical Surveys and Monographs \textbf{81}, American Mathematical Society, Providence, RI, 2000.

\bibitem{T1} 
     	\newblock M. E. Taylor,
     	\newblock \textit{Partial Differential Equations I - III},
     	\newblock Second Edition, Applied Mathematical Sciences \textbf{115} Springer-Verlag, New York, 2011.

\bibitem{Visan}
	\newblock M. Visan,
	\newblock ``The defocusing energy-critical nonlinear {S}chr{\"o}dinger equation in dimensions five and higher'',
	\newblock \textit{PhD Thesis} UCLA (2006).

\bibitem{V2}
	\newblock M. Visan,
	\newblock ``The defocusing energy-critical nonlinear {S}chr{\"o}dinger equation in higher dimensions'',
	\newblock \textit{Duke Mathematical Journal} \textbf{138} (2007) 281 -- 374.

\bibitem{Visan1}
	\newblock M. Visan,
	\newblock ``Global well - posedness and scattering for the defocusing cubic nonlinear {S}chr{\"o}dinger equation in four dimensions",
	\newblock \textit{International Mathematics Research Notices. IMRN} \textbf{5} (2012) 1037 -- 1067.

\bibitem{Yaj}
	\newblock K. Yajima,
	\newblock ``Existence of solutions for {S}chr{\"o}dinger evolution equations'',
	\newblock \textit{Communications in Mathematical Physics} \textbf{110} no. 3 (1987) 415 - 426.

\end{thebibliography}

\end{document}